\documentclass[11pt, a4paper]{article}
\usepackage[utf8]{inputenc}
\usepackage{bbm, bm}
\usepackage[margin = 1in]{geometry}
\usepackage{amsfonts,amsmath,amssymb,amsthm}
\usepackage{mathtools}
\usepackage{mathrsfs}
\usepackage{graphicx,subcaption}
\usepackage{booktabs}
\usepackage{enumitem}
\usepackage{microtype}
\usepackage{tabularx}
\usepackage{hyperref}
\usepackage{multicol}
\hypersetup{
   colorlinks = true,
   citecolor = blue,
   urlcolor = magenta,
   linkcolor = red
}
\usepackage{multirow}
\usepackage[makeroom]{cancel}
\usepackage{tikz}
\usepackage{framed,color}
\definecolor{shadecolor}{rgb}{1,0.8,0.3}
\usepackage{fancybox}
\usepackage{caption}
\usepackage{algorithm}
\usepackage{algorithmic}
\usepackage{aliascnt}
\usepackage{setspace}
\usepackage{longtable}
\allowdisplaybreaks
\usepackage{comment}
\usetikzlibrary{decorations.markings,arrows}

\title{\textbf{Randomized Adaptive Sequential Allocation and the Limits of\\ Bounded Exposure to Inferior Treatments}}
\date{\today}

\newcommand{\Prob}{\mathbb{P}}
\newcommand{\Exp}{\mathbb{E}}
\newcommand{\Var}{\operatorname{Var}}
\newcommand{\ind}{\mathbf{1}}
\newcommand{\F}{\mathcal{F}}
\newcommand{\R}{\mathbb{R}}

\newcommand{\dto}{\xrightarrow{d}}
\newcommand{\pto}{\xrightarrow{p}}

\usepackage{authblk}

\usepackage{natbib}
\setcitestyle{authoryear, open={(},close={)}}

\author[a]{Sampurna Kundu\footnote{\href{mailto:sampurna.kundu58@gmail.com}{sampurna.kundu58@gmail.com}}}
\author[b]{Monitirtha Dey\footnote{\href{mailto:monitirthadey3@gmail.com}{monitirthadey3@gmail.com}}}
\author[a]{Atanu Biswas\footnote{\href{mailto:appubabale@gmail.com}{appubabale@gmail.com}}}

\affil[a]{\footnotesize Indian Statistical Institute, Kolkata, India}
\affil[b]{\footnotesize Indian Statistical Institute, Bangalore, India}

\begin{document}
\maketitle

\theoremstyle{plain}
\newtheorem{theorem}{Theorem}[section]
\newtheorem{lemma}[theorem]{Lemma}
\newtheorem{proposition}[theorem]{Proposition}
\newtheorem*{proposition*}{Proposition}
\newtheorem{corollary}[theorem]{Corollary}
\theoremstyle{definition}
\newtheorem{remark}[theorem]{Remark}
\newtheorem{assumption}[theorem]{Assumption}
\newtheorem{example}[theorem]{Example}

\begin{abstract}
A clinical trial must learn which treatment is superior while sparing its participants from the inferior one.  \citet{KunduJhaBhandari2025} recently showed that a purely deterministic rule, which always assigns the treatment with the larger current sample mean, achieves a strong form of this ethical goal.  The smaller allocation count stays almost surely finite and has all moments bounded uniformly in the sample size.  Deterministic allocation is, however, predictable, and predictability invites selection bias.  We ask whether this bounded-exposure phenomenon survives once the deterministic rule is replaced by a smooth randomized one.  We show that it does not.  Under any fixed non-degenerate symmetric randomization, both treatments are sampled infinitely often, the adaptive sample means are strongly consistent, and the inferior treatment draws a strictly positive limiting allocation proportion.  Consequently, the inferior allocation count grows linearly and its moments are unbounded.  We then study annealed randomization, in which the randomization scale cools with the accumulated information, and obtain summability and divergence criteria that separate finite, sublinear, and linear exposure.  A forced-exploration variant identifies the superior treatment with probability one while keeping the inferior exposure sublinear, and thereby frees the annealed guarantees from any correct-stabilization event.  We also treat normalized multi-treatment rules, ties, indifference zones, unequal variances, general response distributions, and local alternatives.  Finally, we develop the post-trial asymptotic inference that each exposure regime permits.  The results delineate the trade-off between allocation concealment, persistent exploration, and ethical concentration on the superior treatment.
\end{abstract}


\newpage
\section{Introduction}

A clinical trial must serve two objectives at once.  It must learn which treatment is superior, and it must spare its participants from the treatment that is worse.  These two demands pull in opposite directions, since the only way to learn about a treatment is to assign it to patients.  This tension between information and ethics lies at the heart of adaptive sequential design.  It has shaped a vast literature on bandit problems, sequential decision rules, and response-adaptive randomization.

Bandit models furnish one of the earliest frameworks for adaptive experimentation.  The idea of letting the next assignment favour the arm that currently looks better goes back to \citet{Thompson1933}, and the systematic theory of sequential experimentation to \citet{Robbins1952}.  The Bayesian and dynamic-programming treatments of \citet{BerryFristedt1985} study optimal allocation when each observation carries both reward and information.  In a clinical trial, however, pure reward maximization is rarely the goal.  Investigators also need interpretable estimators, valid inference, and protection against operational biases \citep{FriedmanFurbergDeMets2010}.

Response-adaptive randomization is the main line of work that seeks to reconcile these goals.  The play-the-winner idea of \citet{Zelen1969} and the randomized play-the-winner rule of \citet{WeiDurham1978} update the allocation probabilities from the accruing responses.  Subsequent research produced urn designs, optimal-allocation targets, and multi-arm procedures \citep{IvanovaRosenberger2000,RosenbergerEtAl2001,HuRosenberger2006,AtkinsonBiswas2014}.  The doubly adaptive biased coin design of \citet{Eisele1994} aims at a prescribed allocation target, and its limiting proportions were obtained by \citet{HuZhang2004}.  \citet{HuRosenberger2003} compare such procedures through bias, variability, and power, while \citet{BaldiAntogniniGiovagnoli2015} survey the ethics-versus-information trade-off at book length.  A common thread runs through this work.  The assignment stays random.  Randomness guards allocation concealment and curbs selection bias, a concern that goes back to \citet{BlackwellHodges1957}, \citet{Efron1971}, and \citet{Berger2005}.

For normal responses, \citet{BandyopadhyayBiswas2001} proposed an adaptive randomized design whose allocation probability is a smooth function of the current estimated contrast, with room for prognostic adjustment.  A typical two-treatment rule assigns the next unit to treatment $0$ with probability
\[
  \Prob(A_{n+1}=0\mid \F_n)=G(\widehat\theta_{0,n}-\widehat\theta_{1,n}),
\]
where $G$ is symmetric about zero, the probit link $G(x)=\Phi(x/c)$ being a natural choice. Here $\Phi(\cdot)$ is the cdf of standard Normal distribution.  Such a rule favours the treatment that currently looks better, yet it never assigns deterministically.  \citet{BiswasCoad2005} developed related pairwise extensions for several treatments.  These designs are attractive because they are simple, smooth, and genuinely randomized.

A sharper ethical objective was pursued by \citet{Bhandari2007}, \citet{Bhandari2009}, and, most recently, by \citet{KunduJhaBhandari2025}.  For two normal treatments, \citet{KunduJhaBhandari2025} analyse the deterministic rule that always assigns the treatment with the larger current sample mean.  They establish a striking boundedness phenomenon.  When the treatment means are separated, the smaller allocation count is almost surely finite and has all moments bounded uniformly in the sample size.  The same holds, suitably interpreted, for several treatments.  This deterministic rule thus attains the strongest form of bounded exposure one could hope for.

Its very determinism, however, is a practical liability.  Once the data favour one treatment, the other is abandoned, and the next assignment becomes predictable.  Predictable assignments are discouraged in clinical trials.  An investigator who can anticipate the next assignment may, knowingly or not, distort enrolment, timing, or eligibility \citep{BlackwellHodges1957,Berger2005}.  This brings us to the question that motivates the present work.

\medskip
\noindent\emph{Does the bounded-exposure phenomenon of \citet{KunduJhaBhandari2025} survive when the deterministic rule is replaced by a smooth randomized allocation of the form $G(\widehat\theta_{0,n}-\widehat\theta_{1,n})$?}
\medskip

The answer turns out to be negative for every fixed non-degenerate $G$, and the reason is structural.  Suppose $\theta_0>\theta_1$ and the sample means are consistent.  Then the probability of assigning the inferior treatment tends to
\[
  G(\theta_1-\theta_0),
\]
which is strictly positive for any continuous $G$ with full support, the normal cdf among them.  The inferior treatment therefore keeps drawing a positive fraction of all later assignments.  A small scale parameter can shrink this fraction, but it cannot render the inferior count bounded as the trial grows.

The present work studies the randomized counterpart of this bounded-exposure problem in a unified manner.  Its contributions are as follows.  For two treatments, we prove that both allocation counts diverge almost surely, so the adaptive sample means are strongly consistent along the sample paths.  We then establish an almost sure law of large numbers for the allocation proportions.  This shows that the inferior count grows linearly and rules out uniformly bounded moments.  We examine the normalized multi-treatment pairwise rule and find that it, too, grants every treatment a positive limiting proportion.  We isolate the structural reason behind these results, namely that bounded exposure needs a limiting inferior probability equal to zero, which no fixed smooth rule can supply.  We next analyse annealed randomization, in which the scale cools with the accruing information, and classify the inferior exposure as finite, sublinear, or linear through the summability of the tail probabilities.  A forced-exploration safeguard then makes correct identification of the superior treatment an almost sure event, while keeping the inferior exposure sublinear, so that the annealed guarantees no longer rest on a stabilization event.  We also treat ties, indifference zones, unequal variances, general response distributions, and local alternatives.  Finally, we develop the post-trial inference permitted by each regime, and show that the linear, sublinear, and finite regimes lead to qualitatively different inferential conclusions for the inferior treatment.

This paper is structured as follows. We first formally introduce the model and the fixed two-treatment randomized rule, with the relevant notation, in the next section. Section \ref{sec:two} develops the main two-treatment results, namely infinite sampling, consistency, the allocation-proportion limits, the failure of bounded exposure, and the comparison with the deterministic rule. Section \ref{sec:ext-annealed-two} studies annealed randomization, classifies its behaviour through the cooling rate, and presents the forced-exploration safeguard. Section \ref{sec:multi} treats the fixed multi-treatment rule, and Section \ref{sec:ext-annealed-multi} its annealed counterpart. Section \ref{sec:further-extensions} collects further extensions involving ties, indifference zones, unequal variances, general response distributions, and local alternatives. Section \ref{sec:inference} develops the post-trial asymptotic inference under each exposure regime. Section \ref{sec:simulation-study} presents a simulation study and Section \ref{sec:realdata} a real-data illustration. We conclude with a brief discussion in Section \ref{sec:conclusion}. Throughout this work, $\Phi(\cdot)$ denotes the distribution function of the standard normal law and $\ind\{\cdot\}$ an indicator. The symbols $\to$ a.s., $\pto$, and $\dto$ stand for almost sure convergence, convergence in probability, and convergence in distribution, respectively, and $\F_n$ denotes the history of assignments and responses up to stage $n$.

\section{Preliminaries: Model and Notations}
\label{sec:preliminaries}

We first formally introduce the model and fix the notation.  Consider two independent data streams
\[
  X_{0,1},X_{0,2},\ldots \quad \text{and} \quad X_{1,1},X_{1,2},\ldots,
\]
where $X_{j,r}$ has mean $\theta_j$ and variance $\sigma_j^2\in(0,\infty)$ for $j=0,1$.  The normal model $X_{j,r}\sim N(\theta_j,\sigma_j^2)$ is the motivating case.  The main law-of-large-numbers conclusions, however, require only the strong law for each stream.  We assume throughout the two-treatment theory that
\[
  \theta_0>\theta_1,
\]
so treatment 0 is superior and treatment 1 is inferior.

An initial sample of size $m_0\geq 1$ is taken from each treatment.  After $n$ total observations, let $N_{j,n}$ denote the number of observations assigned to treatment $j$, including the initial observations, and let
\[
  \widehat\theta_{j,n}
  =\frac{1}{N_{j,n}}\sum_{r=1}^{N_{j,n}}X_{j,r},\qquad j=0,1,
\]
be the corresponding sample mean.  Let $A_{n+1}\in\{0,1\}$ be the treatment assigned to the next experimental unit, and let $\F_n$ be the $\sigma$-field generated by all assignments and responses observed up to stage $n$.

Let $G$ be a distribution function symmetric about zero, so that
\[
  G(-x)=1-G(x),\qquad x\in\mathbb{R},
\]
at all continuity points of $G$.  The randomized adaptive rule is
\begin{align}
  \Prob(A_{n+1}=0\mid \F_n)&=G(\widehat\theta_{0,n}-\widehat\theta_{1,n}),\label{eq:two-rule-0}\\
  \Prob(A_{n+1}=1\mid \F_n)&=G(\widehat\theta_{1,n}-\widehat\theta_{0,n}).\label{eq:two-rule-1}
\end{align}
The symmetry condition ensures that the two probabilities sum to one.  The common probit choice is
\[
  G(x)=\Phi(x/c),\qquad c>0,
\]
where $c$ controls the degree of randomization.

\begin{assumption}[Non-degenerate fixed randomization]\label{ass:G}
The function $G$ is continuous at $\theta_0-\theta_1$ and $\theta_1-\theta_0$, and
\[
  0<G(x)<1 \quad \text{for every finite } x.
\]
\end{assumption}

Assumption \ref{ass:G} is satisfied by the normal, logistic, and many other smooth symmetric cdf choices.  It excludes the deterministic step rule, which is discussed separately in Section \ref{sec:comparison}.

\section{Main Results for Two Treatments}
\label{sec:two}

The first lemma shows that a fixed non-degenerate randomized rule cannot permanently abandon either treatment.

\begin{lemma}[Infinite sampling of both treatments]\label{lem:both-infinite}
Under Assumption \ref{ass:G},
\[
  N_{0,n}\to\infty \quad \text{and} \quad N_{1,n}\to\infty
  \qquad \text{almost surely.}
\]
\end{lemma}

\begin{proof}
We prove the assertion for treatment 0; the proof for treatment 1 is identical.  Let
\[
  I_{0,n+1}=\ind\{A_{n+1}=0\},
  \qquad p_{0,n}=\Prob(A_{n+1}=0\mid \F_n).
\]
We use the conditional Borel--Cantelli lemma, also known as L\'evy's extension of the Borel--Cantelli lemma; see \citet{Levy1937}, \citet{Freedman1973}, and \citet[Section 12.15]{Williams1991}.  This result does not require independence among the events; the conditioning on the past replaces any independence requirement.  Applied to the adapted events $\{A_{n+1}=0\}$, it implies that the event
\[
  \left\{\sum_{n}p_{0,n}=\infty,\; \sum_n I_{0,n+1}<\infty\right\}
\]
has probability zero.

Suppose, to the contrary, that treatment 0 is sampled only finitely many times.  Then, after some finite random time, $\widehat\theta_{0,n}$ is a fixed finite random variable.  Since the total sample size tends to infinity, treatment 1 must then be sampled infinitely often, and the SLLN gives $\widehat\theta_{1,n}\to\theta_1$ along this event.  Hence
\[
  p_{0,n}=G(\widehat\theta_{0,n}-\widehat\theta_{1,n})
  \longrightarrow G(\widehat\theta_{0,\infty}-\theta_1)>0,
\]
where the last inequality follows from Assumption \ref{ass:G}.  Thus $\sum_n p_{0,n}=\infty$, contradicting conditional Borel--Cantelli.  Therefore treatment 0 is sampled infinitely often almost surely.
\end{proof}

\begin{lemma}[Consistency of adaptive sample means]\label{lem:consistency}
Under the assumptions of Lemma \ref{lem:both-infinite},
\[
  \widehat\theta_{0,n}\to\theta_0,
  \qquad
  \widehat\theta_{1,n}\to\theta_1
  \qquad \text{almost surely.}
\]
\end{lemma}

\begin{proof}
By Lemma \ref{lem:both-infinite}, $N_{j,n}\to\infty$ almost surely for $j=0,1$.  The estimator $\widehat\theta_{j,n}$ is the sample mean of the first $N_{j,n}$ observations from the $j$th independent data stream.  Since the sample mean converges almost surely to $\theta_j$ as its sample size tends to infinity, the result follows by evaluating that convergent sequence at the random times $N_{j,n}$.
\end{proof}

\begin{remark}[No bias from adaptive selection]\label{rem:no-bias}
The step ``evaluate the convergent sequence at the random times $N_{j,n}$'' deserves one word of justification, because $N_{j,n}$ is a data-dependent index.  It is legitimate here because each $\widehat\theta_{j,n}$ is the average of the \emph{first} $N_{j,n}$ entries of a \emph{fixed} i.i.d.\ stream $X_{j,1},X_{j,2},\ldots$: the predictable allocation rule chooses \emph{how many} units arm $j$ receives, not \emph{which} draws from its stream it receives.  Consequently the strong law $m^{-1}\sum_{r=1}^{m}X_{j,r}\to\theta_j$ holds on a single probability-one event for all deterministic $m$ simultaneously, and on that event it continues to hold along any sequence of indices tending to infinity, in particular along $m=N_{j,n}$. Thus, adaptive sampling induces no selection bias in the arm-specific estimators; see \citet{MelfiPage2000} and \citet{LaiWei1982} for general estimation theory after adaptive allocation.  The same remark applies to every later use of the strong law at the random indices $N_{j,n}$.
\end{remark}

\begin{theorem}[Allocation proportions]\label{thm:two-proportions}
Suppose Assumption \ref{ass:G} holds and $\theta_0>\theta_1$.  Then
\[
  \frac{N_{0,n}}{n}\to G(\theta_0-\theta_1),
  \qquad
  \frac{N_{1,n}}{n}\to G(\theta_1-\theta_0)
  \qquad \text{almost surely.}
\]
\end{theorem}

We shall utilize the following result to prove this.

\begin{proposition*}[Martingale SLLN with bounded increments]
Let $(M_n,\mathcal{G}_n)$ be a martingale with $M_0=0$ and increments $D_n=M_n-M_{n-1}$.  If there is a constant $C<\infty$ such that $|D_n|\leq C$ almost surely for every $n$, then
\[
  \frac{M_n}{n}\to0\qquad\text{almost surely.}
\]
This is the bounded-increment form of the martingale SLLN; see, for example, \citet[Chapter 2]{HallHeyde1980} and \citet{Chow1965}.
\end{proposition*}

\begin{proof}
Let $I_{0,t}=\ind\{A_t=0\}$ and write
\[
  N_{0,n}=m_0+\sum_{t=2m_0+1}^{n} I_{0,t},
\]
given the initial $m_0$ observations per treatment.  Define
\[
  p_{0,t-1}=\Prob(A_t=0\mid\F_{t-1})
  =G(\widehat\theta_{0,t-1}-\widehat\theta_{1,t-1}).
\]
The process
\[
  M_n=\sum_{t=2m_0+1}^{n}(I_{0,t}-p_{0,t-1}),
  \qquad n\geq 2m_0+1,
\]
is a martingale with respect to $(\F_n)$.  Indeed, $I_{0,t}$ is $\F_t$-measurable, $p_{0,t-1}$ is $\F_{t-1}$-measurable, and
\[
  \Exp(I_{0,t}-p_{0,t-1}\mid\F_{t-1})
  =\Exp(I_{0,t}\mid\F_{t-1})-p_{0,t-1}=0.
\]
Moreover, since $0\leq I_{0,t}\leq1$ and $0\leq p_{0,t-1}\leq1$, each increment satisfies $|I_{0,t}-p_{0,t-1}|\leq1$.  The martingale SLLN stated above therefore gives $M_n/n\to 0$ almost surely.  By Lemma \ref{lem:consistency},
\[
  p_{0,t-1}\to G(\theta_0-\theta_1) \qquad \text{almost surely},
\]
and therefore Cesaro averaging yields
\[
  \frac{1}{n}\sum_{t=2m_0+1}^{n}p_{0,t-1}
  \to G(\theta_0-\theta_1)\qquad\text{almost surely}.
\]
Combining these facts with $m_0/n\to0$ proves the first assertion.  The second follows either by the same argument or from $N_{1,n}=n-N_{0,n}$ and symmetry of $G$.
\end{proof}

\begin{remark}[Relation to target-allocation procedures]\label{rem:dbcd}
Theorem~\ref{thm:two-proportions} is, in form, a target-allocation statement: the realized proportion $N_{0,n}/n$ converges to the deterministic target $G(\theta_0-\theta_1)$ determined by the rule.  This is the same phenomenon studied for the doubly adaptive biased coin design, whose allocation proportions converge to a prescribed target \citep{Eisele1994,HuZhang2004}.  The distinctive point here is not that a target exists but \emph{where} it sits: for any fixed non-degenerate symmetric $G$ the target assigns the inferior arm the strictly positive mass $G(\theta_1-\theta_0)$, which is exactly what precludes bounded exposure.
\end{remark}

\begin{corollary}[Failure of bounded inferior exposure]\label{cor:two-moments}
Under the assumptions of Theorem \ref{thm:two-proportions},
\[
  \frac{N_{1,n}}{n}\to G(\theta_1-\theta_0)>0
  \qquad \text{almost surely.}
\]
Moreover, for every $r>0$,
\[
  \frac{\Exp\{(N_{1,n})^r\}}{n^r}\to \{G(\theta_1-\theta_0)\}^r.
\]
Thus $\sup_n \Exp\{(N_{1,n})^r\}=\infty$ for every $r>0$.
\end{corollary}

\begin{proof}
The almost sure convergence follows directly from Theorem \ref{thm:two-proportions}.  Since $0\leq N_{1,n}/n\leq1$, dominated convergence gives
\[
  \Exp\left[\left(\frac{N_{1,n}}{n}\right)^r\right]
  \to \{G(\theta_1-\theta_0)\}^r,
\]
which is the desired moment statement.
\end{proof}

\begin{corollary}[The Kundu--Jha--Bhandari count also grows linearly]\label{cor:min-count}
Define the smaller allocation count
\[
  N_{\min,n}=\min\{N_{0,n},N_{1,n}\}.
\]
Suppose Assumption~\ref{ass:G} holds, $G$ is symmetric and \emph{strictly} increasing, and $\theta_0>\theta_1$.  Then
\[
  \frac{N_{\min,n}}{n}\to G(\theta_1-\theta_0)>0
  \qquad \text{almost surely,}
\]
and all positive moments of $N_{\min,n}$ are unbounded in $n$.
\end{corollary}

\begin{proof}
Symmetry gives $G(0)=1/2$.  Since $\theta_0-\theta_1>0$, strict monotonicity yields
\[
  G(\theta_0-\theta_1)>G(0)=\tfrac12>G(\theta_1-\theta_0),
\]
so $G(\theta_1-\theta_0)<1/2<G(\theta_0-\theta_1)$, and Assumption~\ref{ass:G} gives $G(\theta_1-\theta_0)>0$.  By Theorem~\ref{thm:two-proportions}, $N_{0,n}/n\to G(\theta_0-\theta_1)$ and $N_{1,n}/n\to G(\theta_1-\theta_0)$ almost surely, with the second limit strictly smaller; hence $N_{1,n}<N_{0,n}$ for all large $n$ almost surely, so $N_{\min,n}=N_{1,n}$ eventually and $N_{\min,n}/n\to G(\theta_1-\theta_0)$.  The moment statement follows from the dominated-convergence argument used in Corollary~\ref{cor:two-moments}.
\end{proof}

\begin{remark}[Merely monotone $G$]\label{rem:min-monotone}
Strict monotonicity is used only to separate the two limiting proportions.  If $G$ is merely nondecreasing and symmetric (still under Assumption~\ref{ass:G}), then $G(\theta_1-\theta_0)\le 1/2\le G(\theta_0-\theta_1)$ and one always has $\min\{G(\theta_0-\theta_1),G(\theta_1-\theta_0)\}=G(\theta_1-\theta_0)$ because $\theta_1-\theta_0<\theta_0-\theta_1$.  Theorem~\ref{thm:two-proportions} then still gives $N_{\min,n}/n\to G(\theta_1-\theta_0)>0$ almost surely; only the strict ordering $N_{1,n}<N_{0,n}$ may fail to hold eventually (the two counts can be of the same asymptotic size when $G$ is flat at the origin), but the value of the limit is unchanged.
\end{remark}

\begin{remark}[Probit allocation]
For $G(x)=\Phi(x/c)$ and $\Delta=\theta_0-\theta_1>0$, Theorem \ref{thm:two-proportions} gives
\[
  \frac{N_{1,n}}{n}\to \Phi(-\Delta/c).
\]
Choosing a small $c$ can make this limiting fraction very small, but it remains strictly positive for every fixed $c>0$.  Hence the expected number of inferior allocations is asymptotically $n\Phi(-\Delta/c)$, not bounded.
\end{remark}

\subsection{Comparison with the deterministic rule}\label{sec:comparison}

The deterministic allocation rule studied by \citet{KunduJhaBhandari2025} may be viewed formally as the limiting step function
\[
  G_D(x)=
  \begin{cases}
  1, & x>0,\\
  1/2, & x=0,\\
  0, & x<0.
  \end{cases}
\]
This function is symmetric in the sense $G_D(-x)=1-G_D(x)$ for every $x\in\R$ (at $x=0$ both sides equal $1/2$), but it is not a non-degenerate continuous cdf with $0<G_D(x)<1$ at finite negative contrasts.  If $\theta_0>\theta_1$, then
\[
  G_D(\theta_1-\theta_0)=0.
\]
Thus the deterministic design can eventually stop assigning the inferior treatment, while any fixed smooth randomized design assigns it with a positive limiting probability.

This comparison explains why the boundedness theorem of \citet{KunduJhaBhandari2025} does not transfer to the fixed randomized rule.  The bounded-exposure property is not a generic consequence of adaptivity; it is a consequence of a rule whose limiting inferior-treatment allocation probability is zero.  Randomization introduced through a fixed smooth cdf removes that zero and replaces it with a positive exploration rate.

The following proposition summarizes the underlying impossibility principle.

\begin{proposition}[A necessary condition for bounded exposure]\label{prop:necessary}
Consider any two-treatment adaptive design for which the inferior-treatment allocation probability satisfies
\[
  \Prob(A_{n+1}=1\mid\F_n)\to q
  \qquad\text{almost surely}
\]
for some deterministic constant $q>0$ when $\theta_0>\theta_1$.  Then
\[
  \frac{N_{1,n}}{n}\to q
  \qquad\text{almost surely.}
\]
In particular, bounded inferior exposure is impossible unless the limiting inferior allocation probability is zero.
\end{proposition}

\begin{proof}
Write
\[
  N_{1,n}=N_{1,2m_0}+\sum_{t=2m_0+1}^{n} I_{1,t},
\]
where $I_{1,t}=\ind\{A_t=1\}$ and $q_{1,t-1}=\Exp(I_{1,t}\mid\F_{t-1})=\Prob(A_t=1\mid\F_{t-1})$.  The centered process $\sum_{t}(I_{1,t}-q_{1,t-1})$ is a martingale with respect to $(\F_t)$, and its increments satisfy $|I_{1,t}-q_{1,t-1}|\le1$ because both $I_{1,t}$ and $q_{1,t-1}$ lie in $[0,1]$.  Hence the bounded-increment martingale SLLN of the Proposition above applies \emph{unconditionally} (no extra hypothesis is needed), giving
\[
  \frac1n\sum_{t=2m_0+1}^{n}(I_{1,t}-q_{1,t-1})\to0
  \qquad\text{almost surely.}
\]
By assumption $q_{1,t-1}\to q$ almost surely, so Ces\`aro averaging gives $n^{-1}\sum_{t}q_{1,t-1}\to q$, and combining the two displays with $N_{1,2m_0}/n\to0$ proves $N_{1,n}/n\to q$.
\end{proof}

\section{Annealed Randomized Allocation}
\label{sec:ext-annealed-two}

Section~\ref{sec:comparison} shows that fixed smooth randomization cannot reproduce the bounded exposure of the deterministic rule.  The obstruction is not randomness itself.  It is the persistence of a positive inferior probability after the superior treatment has been identified.  A natural remedy is to replace the fixed scale in $G(x/c)$ by a scale that shrinks with the accumulated information.

The idea is to randomize genuinely early in the trial, when the empirical ordering is still uncertain, and to sharpen the rule once the contrast has become reliable.  In this sense the annealed rule interpolates between the fixed randomized design of \citet{BandyopadhyayBiswas2001} and the deterministic rule of \citet{KunduJhaBhandari2025}.  Writing the scale as $c_n$, letting $c_n\downarrow0$ makes a positive contrast increasingly decisive, while still allowing non-degenerate randomization at every finite stage.

This sharpening must be controlled.  If the estimated contrast is eventually positive and bounded away from zero by a constant $d>0$, then the inferior probability is at most the tail term $G(-d/c_n)$.  Finite post-stabilization exposure is thus governed by the summability of these tail probabilities.  A scale that cools too fast, on the other hand, can make the design nearly deterministic before the ordering is trustworthy, and this risks long-run concentration on the wrong treatment.

The sample-size-dependent sharpening mechanism can be formalized as follows. Let $(c_n)$ be a positive deterministic sequence with $c_n\downarrow0$.  The annealed two-treatment rule is
\begin{align}
  \Prob(A_{n+1}=0\mid\F_n)
  &=G\left(\frac{\widehat\theta_{0,n}-\widehat\theta_{1,n}}{c_n}\right),\label{eq:annealed-two-0}\\
  \Prob(A_{n+1}=1\mid\F_n)
  &=G\left(\frac{\widehat\theta_{1,n}-\widehat\theta_{0,n}}{c_n}\right).\label{eq:annealed-two-1}
\end{align}
The deterministic rule of \citet{KunduJhaBhandari2025} may be viewed as the limiting case in which the effective scale has already been reduced to zero.  We note that, relative to the fixed rule of Section~\ref{sec:comparison}, the scale has been moved \emph{out} of $G$: in \eqref{eq:two-rule-0}--\eqref{eq:two-rule-1} the probit choice was written $G(x)=\Phi(x/c)$ with the scale inside, whereas in \eqref{eq:annealed-two-0}--\eqref{eq:annealed-two-1} we take a fixed symmetric cdf $G$ (for the probit rule, $G=\Phi$) and carry the time-varying scale explicitly as the divisor $c_n$.  The two conventions agree when $c_n\equiv c$.

The following theorem gives the basic summability principle.  It is stated on a stabilization event, because annealing can otherwise freeze the design too soon when early estimates mislead.  This qualification is essential, and not a technical artifact.

Throughout this section write
\[
  \widehat\Delta_n=\widehat\theta_{0,n}-\widehat\theta_{1,n}
\]
for the empirical contrast.  For $d>0$ and an integer $T$, define the stabilization event
\[
  \mathcal S(d,T)=\{\widehat\Delta_n\geq d\text{ for all }n\geq T\},
\]
and let
\[
  \mathcal S(d)=\bigcup_{T=1}^{\infty}\mathcal S(d,T).
\]
Thus $\mathcal S(d)$ is the event that the empirical contrast is eventually positive and bounded away from zero by $d$.

\begin{proposition}[Occurrence of stabilization under continued sampling]\label{prop:stabilization-occurs}
Let $\Delta=\theta_0-\theta_1>0$.  Under any adaptive design for which $N_{0,n}\to\infty$, $N_{1,n}\to\infty$ almost surely and the two response streams satisfy the SLLN, one has
\[
  \Prob\{\mathcal S(d)\}=1
  \qquad\text{for every }0<d<\Delta.
\]
In particular, the condition holds for the fixed non-degenerate randomized rule studied in Section~\ref{sec:comparison} and Theorem~\ref{thm:two-proportions}.
\end{proposition}

\begin{proof}
By the SLLN along the random sample sizes, $\widehat\theta_{0,n}\to\theta_0$ and $\widehat\theta_{1,n}\to\theta_1$ almost surely whenever both allocation counts diverge.  Hence $\widehat\Delta_n\to\Delta$ almost surely.  If $0<d<\Delta$, then on this almost sure convergence event there is a finite random time $T$ such that $|\widehat\Delta_n-\Delta|<\Delta-d$ for all $n\geq T$.  Therefore $\widehat\Delta_n>d$ for all $n\geq T$, proving $\mathcal S(d)$.
\end{proof}

\begin{theorem}[Summability criterion after correct stabilization]\label{thm:annealed-summability}
Let $\theta_0>\theta_1$.  Fix $d>0$ and a deterministic integer $T\geq 1$, and consider the annealed rule \eqref{eq:annealed-two-0}--\eqref{eq:annealed-two-1}.  Introduce the stopping time
\[
  \tau=\inf\{n\geq T:\widehat\Delta_n<d\}
  \qquad(\inf\emptyset=\infty),
\]
so that $\mathcal S(d,T)=\{\tau=\infty\}$, and write $\widetilde I_{n+1}=\ind\{A_{n+1}=1\}\,\ind\{\tau>n\}$ for the inferior-allocation indicators counted only while the contrast stays above $d$.  If
\begin{equation}\label{eq:annealed-summable-condition}
  \sum_{n=1}^{\infty}G(-d/c_n)<\infty,
\end{equation}
then for every $t>0$,
\[
  \Exp\left[\exp\left\{t\sum_{n=T}^{\infty}\widetilde I_{n+1}\right\}\,\middle|\,\F_T\right]
  \leq
  \exp\left\{(e^t-1)\sum_{n=T}^{\infty}G(-d/c_n)\right\}<\infty
  \qquad\text{a.s.}
\]
Consequently $\sum_{n\geq T}\widetilde I_{n+1}<\infty$ almost surely with all moments finite, and on the stabilization event $\mathcal S(d,T)$ the number of inferior-treatment allocations after time $T$ equals $\sum_{n\geq T}\widetilde I_{n+1}$ and is therefore finite almost surely.  Taking the union over $T$, the post-stabilization inferior count is finite almost surely on $\mathcal S(d)=\bigcup_{T\geq1}\mathcal S(d,T)$.
\end{theorem}

\begin{proof}
\emph{Step 1: $\tau$ is a stopping time and the per-step bound is adapted.}  For $n\geq T$,
\[
  \{\tau>n\}=\{\widehat\Delta_m\geq d\text{ for all }T\leq m\leq n\}\in\F_n ,
\]
because each $\widehat\Delta_m$ is $\F_m$-measurable; hence $\tau$ is an $(\F_n)$-stopping time and $\{\tau=\infty\}=\mathcal S(d,T)$.  This is the device that repairs the measurability gap: the event $\mathcal S(d,T)$ depends on the whole future and is not $\F_T$-measurable, but the truncated event $\{\tau>n\}$ is $\F_n$-measurable, so it may be used inside conditional expectations given $\F_n$.

\emph{Step 2: conditional mean of the truncated indicators.}  Since $\widetilde I_{n+1}=\ind\{A_{n+1}=1\}\ind\{\tau>n\}$ and $\ind\{\tau>n\}$ is $\F_n$-measurable,
\[
  \Exp(\widetilde I_{n+1}\mid\F_n)
  =\ind\{\tau>n\}\,G\!\left(-\frac{\widehat\Delta_n}{c_n}\right)
  \leq \ind\{\tau>n\}\,G(-d/c_n)
  \leq G(-d/c_n)=:q_n ,
\]
where the first inequality holds because on $\{\tau>n\}$ we have $\widehat\Delta_n\geq d$ and $G$ is nondecreasing, so $G(-\widehat\Delta_n/c_n)\leq G(-d/c_n)$.

\emph{Step 3: exponential-moment telescoping.}  Fix $t>0$.  As $\widetilde I_{M+1}\in\{0,1\}$, $e^{t\widetilde I_{M+1}}=1+(e^t-1)\widetilde I_{M+1}$, hence by Step 2,
\[
  \Exp\!\left(e^{t\widetilde I_{M+1}}\mid\F_M\right)
  =1+(e^t-1)\Exp(\widetilde I_{M+1}\mid\F_M)
  \leq 1+(e^t-1)q_M .
\]
Therefore, using the tower property and that $1+(e^t-1)q_M$ is deterministic,
\[
\begin{aligned}
  \Exp\!\left[\exp\!\Big(t\sum_{n=T}^{M}\widetilde I_{n+1}\Big)\,\middle|\,\F_T\right]
  &=\Exp\!\left[\exp\!\Big(t\sum_{n=T}^{M-1}\widetilde I_{n+1}\Big)\,
      \Exp\!\big(e^{t\widetilde I_{M+1}}\mid\F_M\big)\,\middle|\,\F_T\right]\\
  &\leq \big(1+(e^t-1)q_M\big)\,
      \Exp\!\left[\exp\!\Big(t\sum_{n=T}^{M-1}\widetilde I_{n+1}\Big)\,\middle|\,\F_T\right].
\end{aligned}
\]
Iterating and using $1+x\leq e^x$ gives
\[
  \Exp\!\left[\exp\!\Big(t\sum_{n=T}^{M}\widetilde I_{n+1}\Big)\,\middle|\,\F_T\right]
  \leq \prod_{n=T}^{M}\{1+(e^t-1)q_n\}
  \leq \exp\!\Big\{(e^t-1)\sum_{n=T}^{M}q_n\Big\}.
\]
Letting $M\to\infty$ and invoking the monotone convergence theorem on the left and \eqref{eq:annealed-summable-condition} on the right proves the displayed exponential-moment bound.

\emph{Step 4: conclusion.}  A nonnegative random variable with a finite exponential moment is finite almost surely and has finite moments of every order; applied to $\sum_{n\geq T}\widetilde I_{n+1}$ this gives the moment claim.  Finally, on $\mathcal S(d,T)=\{\tau=\infty\}$ we have $\ind\{\tau>n\}=1$ for all $n\geq T$, so $\widetilde I_{n+1}=\ind\{A_{n+1}=1\}$ there and $\sum_{n\geq T}\widetilde I_{n+1}$ is exactly the post-$T$ inferior count.  Since $\mathcal S(d)$ is the countable union of the $\mathcal S(d,T)$, the post-stabilization inferior count is finite almost surely on $\mathcal S(d)$.
\end{proof}

\begin{theorem}[Divergence criterion after stable separation]\label{thm:annealed-divergence}
Assume the setup of Theorem \ref{thm:annealed-summability}.  Suppose, in addition, that for some finite $D>d$ and some almost surely finite time $T$,
\[
  d\leq \widehat\Delta_n\leq D\qquad\text{for all }n\geq T.
\]
If
\[
  \sum_{n=1}^{\infty}G(-D/c_n)=\infty,
\]
then the number of inferior-treatment allocations after $T$ is infinite almost surely.
\end{theorem}

\begin{proof}
For $n\geq T$,
\[
  \Prob(A_{n+1}=1\mid\F_n)
  =G(-\widehat\Delta_n/c_n)
  \geq G(-D/c_n).
\]
The lower bound has divergent series.  The conditional Borel--Cantelli lemma therefore implies that $A_{n+1}=1$ occurs infinitely often almost surely on the stated stable-separation event.
\end{proof}

\begin{remark}[Why stabilization has to be mentioned]
If $c_n$ decreases very rapidly, the rule can become nearly deterministic before the empirical ordering is reliable.  On paths where the early empirical ordering is wrong, the design may allocate almost all later patients to the wrong treatment.  Therefore annealing needs either a sufficiently reliable pilot stage, a delayed cooling schedule, or a forced-exploration safeguard.  This point is analogous to the usual exploration--exploitation issue in bandit problems; see, for example, \citet{LaiRobbins1985}, \citet{AuerCesaBianchiFischer2002}, \citet{BubeckCesaBianchi2012}, and \citet{LattimoreSzepesvari2020}.
\end{remark}

\subsection{An unconditional guarantee via forced exploration}\label{sec:forced-exploration}

Theorems~\ref{thm:annealed-summability} and \ref{thm:annealed-divergence} are stated \emph{on} the stabilization event $\mathcal S(d,T)$.  This is appropriate, but it leaves open the question raised in the preceding remark: for the pure annealed rule one does not control $\Prob\{\mathcal S(d)\}$, because fast cooling can lock the design onto the wrong arm.  We now show that a small amount of forced exploration removes this gap entirely: it forces stabilization on the \emph{correct} ordering with probability one, so the conditional conclusions above become unconditional, while inferior exposure stays sublinear.

Let $(\varepsilon_n)\subset(0,1)$ be a deterministic sequence with
\[
  \varepsilon_n\downarrow0,\qquad
  \sum_{n}\varepsilon_n=\infty,\qquad
  E_N:=\sum_{n=1}^{N}\varepsilon_n=o(N).
\]
Let $(U_n)$ be i.i.d.\ $\mathrm{Uniform}(0,1)$ random variables, independent of the two response streams, and enlarge the filtration so that $\F_n$ contains $U_1,\ldots,U_n$ as well.  The \emph{guarded annealed rule} assigns the $(n{+}1)$st unit as follows: if $U_{n+1}\leq\varepsilon_{n+1}$ (a forced-exploration step) it assigns treatment $0$ or $1$ with probability $1/2$ each; otherwise it uses the annealed probit allocation $\Phi(\widehat\Delta_n/c_n)$ for treatment $0$.  Equivalently,
\begin{equation}\label{eq:guarded-rule}
  \Prob(A_{n+1}=1\mid\F_n)
  =\tfrac12\,\varepsilon_{n+1}+(1-\varepsilon_{n+1})\,\Phi\!\left(-\frac{\widehat\Delta_n}{c_n}\right),
\end{equation}
and symmetrically for treatment $0$.

\begin{theorem}[Unconditional consistency with sublinear inferior exposure]\label{thm:forced-exploration}
Assume $\theta_0>\theta_1$, finite variances, and the guarded annealed rule \eqref{eq:guarded-rule} with a cooling sequence $c_n\downarrow0$ satisfying $\sum_n\Phi(-d/c_n)<\infty$ for some $0<d<\Delta:=\theta_0-\theta_1$.  Then, with no conditioning on any stabilization event:
\begin{enumerate}[label=\textnormal{(\roman*)}]
\item \textnormal{(Unconditional correct stabilization.)}  Both treatments are sampled infinitely often, $\widehat\theta_{j,n}\to\theta_j$ almost surely for $j=0,1$, and $\widehat\Delta_n\to\Delta>0$ almost surely.  In particular $\Prob\{\mathcal S(d)\}=1$ for every $0<d<\Delta$, so the conclusion of Theorem~\ref{thm:annealed-summability} holds almost surely rather than merely on a stabilization event.
\item \textnormal{(Sublinear exposure.)}  The inferior-allocation count satisfies
\[
  N_{1,N}\sim \tfrac12 E_N\qquad\text{almost surely,}
\]
hence $N_{1,N}\to\infty$ and $N_{1,N}/N\to0$ almost surely.  The superior treatment receives $N_{0,N}=N-N_{1,N}\sim N$.
\end{enumerate}
\end{theorem}

\begin{proof}
\emph{(i)}  By \eqref{eq:guarded-rule} and its symmetric counterpart, both conditional allocation probabilities are bounded below by the forced-exploration floor:
\[
  \Prob(A_{n+1}=0\mid\F_n)\geq\tfrac12\varepsilon_{n+1},
  \qquad
  \Prob(A_{n+1}=1\mid\F_n)\geq\tfrac12\varepsilon_{n+1}.
\]
Since $\sum_n\tfrac12\varepsilon_{n+1}=\infty$, the conditional Borel--Cantelli lemma (as in Lemma~\ref{lem:both-infinite}) applied to each arm gives $N_{0,n}\to\infty$ and $N_{1,n}\to\infty$ almost surely; the floor is deterministic, so no appeal to the empirical contrast is needed and the argument is unconditional.  Consistency $\widehat\theta_{j,n}\to\theta_j$ then follows as in Lemma~\ref{lem:consistency} (using Remark~\ref{rem:no-bias} for the random indices), whence $\widehat\Delta_n\to\Delta>0$.  By Proposition~\ref{prop:stabilization-occurs}, $\Prob\{\mathcal S(d)\}=1$ for every $0<d<\Delta$, so the post-stabilization finiteness of Theorem~\ref{thm:annealed-summability} now holds on a probability-one event.

\emph{(ii)}  Let $I_{n+1}=\ind\{A_{n+1}=1\}$ and $q_n=\Prob(A_{n+1}=1\mid\F_n)$ from \eqref{eq:guarded-rule}.  Write
\[
  q_n=\tfrac12\varepsilon_{n+1}+(1-\varepsilon_{n+1})\Phi(-\widehat\Delta_n/c_n),
  \qquad
  Q_N=\sum_{n=2m_0}^{N-1}q_n .
\]
By part (i) there is an almost surely finite $T$ with $\widehat\Delta_n\geq d$ for $n\geq T$; for such $n$, monotonicity of $\Phi$ gives $\Phi(-\widehat\Delta_n/c_n)\leq\Phi(-d/c_n)$, and $\sum_n\Phi(-d/c_n)<\infty$ by hypothesis.  Hence the annealed contribution $\sum_n(1-\varepsilon_{n+1})\Phi(-\widehat\Delta_n/c_n)$ converges almost surely, and
\[
  Q_N=\tfrac12\sum_{n=2m_0}^{N-1}\varepsilon_{n+1}+O(1)
      =\tfrac12 E_N+O(1)\qquad\text{a.s.},
\]
so $Q_N\sim\tfrac12 E_N\to\infty$.  The centered process $M_N=\sum_{n=2m_0}^{N-1}(I_{n+1}-q_n)$ is a martingale difference sum with $|I_{n+1}-q_n|\le1$ and conditional variance $\Exp\{(I_{n+1}-q_n)^2\mid\F_n\}=q_n(1-q_n)\le q_n$, so $\sum_n\Var(\cdot\mid\F_n)\le Q_N$.  Because $Q_N\uparrow\infty$, the Chow martingale strong law \citep{Chow1965} (equivalently the square-integrable martingale convergence theorem applied to $\sum (I_{n+1}-q_n)/Q_n$ followed by Kronecker's lemma, exactly as in Theorem~\ref{thm:cooling-classification}(iii)) gives $M_N/Q_N\to0$ almost surely; the required normalization $\sum_N q_N/Q_N^2<\infty$ holds automatically since $Q_N\uparrow\infty$.  Therefore
\[
  N_{1,N}=N_{1,2m_0}+\sum_{n=2m_0}^{N-1}I_{n+1}
         =Q_N+M_N\sim Q_N\sim\tfrac12 E_N
  \qquad\text{a.s.}
\]
Since $E_N=o(N)$, we get $N_{1,N}/N\to0$; and $E_N\to\infty$ gives $N_{1,N}\to\infty$.  The complementary statement for $N_{0,N}$ is immediate from $N_{0,N}=N-N_{1,N}$.
\end{proof}

\begin{remark}[Reading of the guarantee, and a pilot-stage alternative]\label{rem:forced-reading}
Taking, for instance, $\varepsilon_n=n^{-1/2}$ gives $E_N\sim2\sqrt N$ and hence $N_{1,N}\sim\sqrt N$: the inferior arm is visited infinitely often (which, as Section~\ref{sec:inference} shows, is exactly what makes inference about $\theta_1$ possible) but only at a vanishing rate, and the design identifies the superior treatment with probability one.  Thus forced exploration converts the conditional finite-/sublinear-exposure statements of Section~\ref{sec:ext-cooling-classification} into an unconditional design with provable correct selection.  A two-phase alternative achieves the same end: run a deterministic pilot of length $m(N)\to\infty$ with $m(N)/N\to0$ to fix the ordering---Lemma~\ref{lem:general-stabilization} bounds the probability of a wrong ordering by an exponentially small term---and anneal thereafter; the pilot contributes $o(N)$ inferior observations while guaranteeing $\Prob(\text{correct selection})\to1$.
\end{remark}

\subsection{Classification according to the cooling rate for \texorpdfstring{$G(\cdot)=\Phi(\cdot)$}{G(.)=Phi(.)}}\label{sec:ext-cooling-classification}

The term \emph{cooling rate} refers to the speed at which the scale $c_n$ decreases to zero.  A faster cooling rate makes the allocation rule closer to the deterministic step rule sooner, whereas a slower cooling rate preserves randomization for a longer period.  A \emph{cooling regime} is an asymptotic class of schedules $c_n$ characterized by whether the associated inferior-tail probabilities are summable, divergent but Cesaro negligible, or bounded away from zero.  The previous section shows that the behavior of the annealed rule is governed by the series $\sum_n G(-d/c_n)$.  For the probit rule $G(x)=\Phi(x)$, this condition can be made explicit.

\begin{theorem}[Three cooling regimes for probit annealing]\label{thm:cooling-classification}
Assume $G=\Phi$ and suppose that, after some almost surely finite time $T$,
\[
  d\leq \widehat\Delta_n\leq D
  \qquad(n\geq T)
\]
for constants $0<d\leq D<\infty$.  Let $N_{1,N}$ denote the inferior-treatment allocation count up to time $N$.

\begin{enumerate}[label=(\roman*)]
\item If $c_n\equiv c>0$, then $N_{1,N}/N$ converges to a positive constant whenever $\widehat\Delta_n$ converges.  In particular, fixed randomization gives linear inferior exposure.

\item If $c_n\downarrow0$ and $\sum_n \Phi(-d/c_n)<\infty$, then $N_{1,\infty}<\infty$ almost surely on the stable-separation event.

\item Put
\[
  a_n=\Phi(-d/c_n),\qquad b_n=\Phi(-D/c_n),
\]
and suppose $c_n\downarrow0$, $\sum_n b_n=\infty$, and $N^{-1}\sum_{n=1}^{N}a_n\to0$.  Then inferior allocations are infinite but have vanishing proportion:
\[
  N_{1,N}\to\infty\quad\text{a.s.},
  \qquad
  \frac{N_{1,N}}{N}\to0\quad\text{a.s.}
\]
Moreover, if $B_N=\sum_{n=T}^{N-1}b_n$ and
\[
  \sum_{N=1}^{\infty}\frac{a_N}{B_N^2}<\infty,
\]
then the martingale fluctuation around the accumulated conditional inferior-allocation probability is $o(B_N)$ almost surely.
\end{enumerate}
\end{theorem}

\begin{proof}
Let
\[
  I_{n+1}=\ind\{A_{n+1}=1\},
  \qquad
  q_n=\Prob(A_{n+1}=1\mid\F_n)=\Phi(-\widehat\Delta_n/c_n),
  \qquad n\geq T.
\]
On the event $d\leq\widehat\Delta_n\leq D$ for all $n\geq T$, monotonicity of $\Phi$ gives
\[
  b_n=\Phi(-D/c_n)\leq q_n\leq \Phi(-d/c_n)=a_n.
\]

If $c_n\equiv c>0$ and $\widehat\Delta_n\to\Delta>0$, then $q_n\to\Phi(-\Delta/c)>0$.  The same bounded-increment martingale SLLN used in Theorem~\ref{thm:two-proportions} gives
\[
  \frac{N_{1,N}}{N}
  =\frac{1}{N}\sum_{n<T} I_{n+1}
   +\frac{1}{N}\sum_{n=T}^{N-1}q_n
   +o(1)
  \to \Phi(-\Delta/c)>0,
\]
which proves part (i).

Part (ii) follows from Theorem~\ref{thm:annealed-summability}, since $q_n\leq a_n$ and $\sum_n a_n<\infty$.

For part (iii), $q_n\geq b_n$ and $\sum_n b_n=\infty$, so $\sum_n q_n=\infty$.  L\'evy's conditional Borel--Cantelli lemma therefore implies that $I_{n+1}=1$ infinitely often, and hence $N_{1,N}\to\infty$ almost surely.  Also,
\[
  \frac{1}{N}\sum_{n=T}^{N-1}q_n
  \leq \frac{1}{N}\sum_{n=T}^{N-1}a_n\to0.
\]
The centered process
\[
  M_N=\sum_{n=T}^{N-1}(I_{n+1}-q_n)
\]
is a martingale with bounded increments, because $I_{n+1}$ is $\F_{n+1}$-measurable, $q_n$ is $\F_n$-measurable, $\Exp(I_{n+1}-q_n\mid\F_n)=0$, and $|I_{n+1}-q_n|\leq1$.  The martingale SLLN gives $M_N/N\to0$ almost surely.  Therefore
\[
  \frac{N_{1,N}}{N}
  =o(1)+\frac{1}{N}\sum_{n=T}^{N-1}q_n+\frac{M_N}{N}\to0
  \qquad\text{a.s.}
\]
Finally, assume $\sum_N a_N/B_N^2<\infty$.  Put
\[
  X_N=I_{N+1}-q_N,\qquad N\geq T,
\]
so that $(X_N,\F_{N+1})$ is a martingale-difference sequence and
\[
  \Exp(X_N^2\mid\F_N)
  =q_N(1-q_N)\leq q_N\leq a_N .
\]
Therefore
\[
  \sum_{N=T}^{\infty}
  \Exp\left(\left.\frac{X_N^2}{B_N^2}\right|\F_N\right)
  \leq
  \sum_{N=T}^{\infty}\frac{a_N}{B_N^2}<\infty .
\]
By the martingale convergence theorem applied to the square-integrable martingale
\[
  S_m=\sum_{N=T}^{m}\frac{X_N}{B_N},
\]
the series $\sum_{N=T}^{\infty}X_N/B_N$ converges almost surely.  Since $B_N$ is nondecreasing and $B_N\to\infty$, Kronecker's lemma yields
\[
  \frac{1}{B_N}\sum_{n=T}^{N-1}X_n\to0
  \qquad\text{almost surely.}
\]
But $\sum_{n=T}^{N-1}X_n=M_N$.  Hence $M_N/B_N\to0$ almost surely, which proves the final assertion.
\end{proof}

\begin{corollary}[The boundary for $c_n=\gamma/\sqrt{\log n}$]\label{cor:probit-log-boundary}
Let $G=\Phi$, let the stable contrast satisfy $\widehat\Delta_n\to\Delta>0$, and take
\[
  c_n=\frac{\gamma}{\sqrt{\log n}},\qquad \gamma>0.
\]
Then
\[
  \Phi(-\Delta/c_n)
  \sim
  \frac{\gamma}{\Delta\sqrt{2\pi\log n}}
  n^{-\Delta^2/(2\gamma^2)}.
\]
Thus the post-stabilization inferior allocation series is summable if
\[
  \gamma<\frac{\Delta}{\sqrt2},
\]
and divergent if
\[
  \gamma\geq\frac{\Delta}{\sqrt2}.
\]
At the boundary $\gamma=\Delta/\sqrt2$, the expected number of post-stabilization inferior allocations grows on the order of $\sqrt{\log N}$.
\end{corollary}

\begin{proof}
The displayed asymptotic equivalence is Mills' ratio for the standard normal tail, $\Phi(-x)\sim\phi(x)/x$ as $x\to\infty$, applied with $x=\Delta/c_n=\Delta\sqrt{\log n}/\gamma$.

\emph{Replacing $\widehat\Delta_n$ by $\Delta$ is legitimate.}  The per-step inferior probability is $q_n=\Phi(-\widehat\Delta_n/c_n)$, not $\Phi(-\Delta/c_n)$, and because the probit tail is exponentially sensitive to its argument while $c_n\downarrow0$, this substitution requires justification.  Write
\[
  \frac{\Phi(-\widehat\Delta_n/c_n)}{\Phi(-\Delta/c_n)}
  =\exp\!\left\{\frac{\Delta^2-\widehat\Delta_n^2}{2c_n^2}\bigl(1+o(1)\bigr)\right\},
  \qquad
  \frac{\Delta^2-\widehat\Delta_n^2}{2c_n^2}
  =\frac{(\Delta-\widehat\Delta_n)(\Delta+\widehat\Delta_n)}{2}\cdot\frac{\log n}{\gamma^2},
\]
where the first equality uses Mills' ratio for both tails.  By the law of the iterated logarithm applied to the superior- and inferior-arm sample means at their (linearly growing) sample sizes, $|\widehat\Delta_n-\Delta|=O\!\big(\sqrt{\log\log n/n}\big)$ almost surely, so
\[
  \left|\frac{\Delta^2-\widehat\Delta_n^2}{2c_n^2}\right|
  =O\!\left(\frac{\log n}{\gamma^2}\sqrt{\frac{\log\log n}{n}}\right)
  \longrightarrow 0
  \qquad\text{a.s.}
\]
Hence $q_n/\Phi(-\Delta/c_n)\to1$ almost surely, and by the limit-comparison test the series $\sum_n q_n$ and $\sum_n\Phi(-\Delta/c_n)$ converge or diverge together.  It therefore suffices to analyze $\sum_n\Phi(-\Delta/c_n)$.

\emph{The series test.}  By the displayed equivalence, $\sum_n\Phi(-\Delta/c_n)$ has the same behavior as
\[
  \sum_n \frac{n^{-\Delta^2/(2\gamma^2)}}{\sqrt{\log n}}.
\]
This series converges exactly when $\Delta^2/(2\gamma^2)>1$, i.e.\ $\gamma<\Delta/\sqrt2$, and diverges when $\gamma\geq\Delta/\sqrt2$.  At equality $\Delta^2/(2\gamma^2)=1$ it is comparable to $\sum_n\{n\sqrt{\log n}\}^{-1}$, whose partial sums are of order $\sqrt{\log N}$ (compare with $\int^N (x\sqrt{\log x})^{-1}\,dx=2\sqrt{\log N}$).  Combined with $\sum_n q_n\sim\sum_n\Phi(-\Delta/c_n)$, this gives the stated summability dichotomy and the $\sqrt{\log N}$ growth of the expected post-stabilization inferior count at the boundary.
\end{proof}

\section{Multi-treatment Randomized Allocation}
\label{sec:multi}

We now pass from two treatments to several.  The randomized rule extends naturally through pairwise comparisons.  The extension, however, needs care, since a naive pairwise score is not a probability vector once there are more than two treatments.

Suppose there are $k\geq2$ treatments with independent data streams
\[
  X_{j,1},X_{j,2},\ldots,\qquad j=1,\ldots,k,
\]
where $X_{j,r}$ has mean $\theta_j$ and finite positive variance.  Assume, for clarity of exposition, that
\[
  \theta_1>\theta_2>\cdots>\theta_k.
\]
Let $N_{j,n}$ and $\widehat\theta_{j,n}$ denote the allocation count and current sample mean for treatment $j$.

A natural pairwise score for treatment $j$ is
\[
  W_{j,n}=\sum_{i\neq j} G(\widehat\theta_{j,n}-\widehat\theta_{i,n}).
\]
The expression $W_{j,n}/(k-1)$ is the average pairwise probability that treatment $j$ beats a randomly selected opponent.  However, it is not a probability vector over $j=1,\ldots,k$ when $k>2$.  Indeed, by pairwise symmetry,
\[
  \sum_{j=1}^k W_{j,n}
  =\sum_{1\leq i<j\leq k}\left\{G(\widehat\theta_{j,n}-\widehat\theta_{i,n})+G(\widehat\theta_{i,n}-\widehat\theta_{j,n})\right\}
  =\binom{k}{2}.
\]
Therefore
\[
  \sum_{j=1}^k \frac{W_{j,n}}{k-1}=\frac{k}{2},
\]
which equals one only when $k=2$.  To obtain a genuine allocation probability vector, the pairwise scores must be normalized as
\begin{equation}\label{eq:multi-rule}
  \pi_{j,n}
  =\frac{W_{j,n}}{\binom{k}{2}}
  =\frac{2}{k(k-1)}\sum_{i\neq j}G(\widehat\theta_{j,n}-\widehat\theta_{i,n}),
  \qquad j=1,\ldots,k.
\end{equation}
This normalized rule reduces exactly to the two-treatment rule when $k=2$.

\subsection{Limit theory for the normalized multi-treatment rule}

Let $A_{n+1}\in\{1,\ldots,k\}$ be the treatment assigned at the next stage, with
\[
  \Prob(A_{n+1}=j\mid\F_n)=\pi_{j,n},
\]
where $\pi_{j,n}$ is given by \eqref{eq:multi-rule}.  Assume initially that each treatment has at least one observation.

\begin{lemma}[Infinite sampling in the multi-treatment design]\label{lem:multi-infinite}
Under Assumption \ref{ass:G}, for each $j=1,\ldots,k$,
\[
  N_{j,n}\to\infty\qquad\text{almost surely.}
\]
Consequently,
\[
  \widehat\theta_{j,n}\to\theta_j,
  \qquad j=1,\ldots,k,
  \qquad\text{almost surely.}
\]
\end{lemma}

\begin{proof}
Fix $j$.  Let $I_{j,n+1}=\ind\{A_{n+1}=j\}$ and let $\pi_{j,n}=\Exp(I_{j,n+1}\mid\F_n)$.  If treatment $j$ were sampled only finitely many times, then $\widehat\theta_{j,n}$ would eventually be a fixed finite random variable.  For every other treatment $i$, either $N_{i,n}$ is finite, in which case $\widehat\theta_{i,n}$ is eventually fixed and finite, or $N_{i,n}\to\infty$, in which case the SLLN gives $\widehat\theta_{i,n}\to\theta_i$.  In all cases, each pairwise term $G(\widehat\theta_{j,n}-\widehat\theta_{i,n})$ has a positive finite limit or positive finite limit inferior.  Hence $\pi_{j,n}$ has a positive limit inferior, so $\sum_n \pi_{j,n}=\infty$ on the event that $N_{j,n}$ is finite.  Conditional Borel--Cantelli then rules out finite sampling of treatment $j$.  Since $j$ was arbitrary, all treatment counts diverge almost surely.  Consistency follows by applying the SLLN to each independent stream at the random sample sizes $N_{j,n}$.
\end{proof}

Define the limiting allocation probabilities
\begin{equation}\label{eq:pi-star}
  \pi_j^*
  =\frac{2}{k(k-1)}\sum_{i\neq j}G(\theta_j-\theta_i),
  \qquad j=1,\ldots,k.
\end{equation}

\begin{theorem}[Multi-treatment allocation proportions]\label{thm:multi-proportions}
Under the normalized rule \eqref{eq:multi-rule} and Assumption \ref{ass:G},
\[
  \frac{N_{j,n}}{n}\to \pi_j^*,
  \qquad j=1,\ldots,k,
  \qquad\text{almost surely.}
\]
Each $\pi_j^*$ is strictly positive.  Hence every treatment, including every inferior treatment, receives a positive limiting fraction of the total sample.
\end{theorem}

\begin{proof}
The proof is the same martingale-averaging argument used in Theorem \ref{thm:two-proportions}.  For each fixed $j$,
\[
  N_{j,n}=N_{j,km_0}+\sum_{t=km_0+1}^{n} I_{j,t},
  \qquad
  \Exp(I_{j,t}\mid\F_{t-1})=\pi_{j,t-1}.
\]
The centered sum is a bounded-increment martingale and is $o(n)$ almost surely.  Lemma \ref{lem:multi-infinite} implies $\widehat\theta_{j,n}\to\theta_j$ for every $j$, so $\pi_{j,n}\to\pi_j^*$ almost surely.  Cesaro averaging completes the proof.  Strict positivity follows from $0<G(x)<1$ for finite $x$.
\end{proof}

\begin{proposition}[Ordering of limiting proportions]\label{prop:ordering}
If $G$ is strictly increasing and $\theta_1>\theta_2>\cdots>\theta_k$, then
\[
  \pi_1^*>\pi_2^*>\cdots>\pi_k^*.
\]
Thus the best treatment has the largest limiting allocation proportion, but every inferior treatment has a positive limiting allocation proportion.
\end{proposition}

\begin{proof}
For $a<b$, compare the pairwise scores corresponding to treatments $a$ and $b$.  Since $\theta_a>\theta_b$,
\[
  G(\theta_a-\theta_b)>G(\theta_b-\theta_a).
\]
For any $\ell\notin\{a,b\}$,
\[
  \theta_a-\theta_\ell>\theta_b-\theta_\ell,
\]
so strict monotonicity gives
\[
  G(\theta_a-\theta_\ell)>G(\theta_b-\theta_\ell).
\]
Adding these inequalities over all relevant pairwise terms yields $\pi_a^*>\pi_b^*$.  Applying this to adjacent pairs proves the assertion.
\end{proof}

Let
\[
  V_n=\sum_{j=2}^{k}N_{j,n}
\]
be the total number of allocations to inferior treatments, and let
\[
  N^{(2)}_n \coloneqq \text{the second-largest value among } N_{1,n},\ldots,N_{k,n}.
\]
The next corollary gives the multi-treatment analogue of Corollary \ref{cor:two-moments}.

\begin{corollary}[Linear suboptimal allocation in the multi-treatment design]\label{cor:multi-linear}
Under the assumptions of Theorem \ref{thm:multi-proportions},
\[
  \frac{V_n}{n}\to 1-\pi_1^*>0
  \qquad\text{almost surely.}
\]
If, in addition, $G$ is strictly increasing and $\theta_1>\theta_2>\cdots>\theta_k$, then
\[
  \frac{N^{(2)}_n}{n}\to \pi_2^*>0
  \qquad\text{almost surely.}
\]
For every $r>0$,
\[
  \frac{\Exp\{(V_n)^r\}}{n^r}\to (1-\pi_1^*)^r,
  \qquad
  \frac{\Exp\{(N^{(2)}_n)^r\}}{n^r}\to (\pi_2^*)^r,
\]
where the second convergence assumes strict ordering.  Hence neither the total inferior-treatment count nor the second-largest allocation count has uniformly bounded moments.
\end{corollary}

\begin{proof}
The almost sure limits follow from Theorem \ref{thm:multi-proportions} and Proposition \ref{prop:ordering}.  The moment limits follow by dominated convergence after division by $n^r$, because $0\leq V_n/n\leq1$ and $0\leq N^{(2)}_n/n\leq1$.
\end{proof}

\section{Annealed Multi-treatment Allocation}
\label{sec:ext-annealed-multi}

There are two distinct annealed multi-treatment questions.  The first asks what happens if the corrected pairwise-sum rule is annealed directly.  The second asks how to modify the multi-treatment rule so that the total allocation to all inferior treatments can be finite or sublinear.

\subsection{Direct annealing of the pairwise-sum rule}

Consider
\[
  W_{j,n}(c_n)=\sum_{i\neq j}G\left(\frac{\widehat\theta_{j,n}-\widehat\theta_{i,n}}{c_n}\right),
  \qquad
  \pi_{j,n}(c_n)=\frac{W_{j,n}(c_n)}{\binom{k}{2}}.
\]
This remains a probability vector because pairwise symmetry still gives $\sum_j W_{j,n}(c_n)=\binom{k}{2}$.

\begin{proposition}[Rank limit of direct pairwise annealing]\label{prop:rank-limit-pairwise}
Suppose $\theta_1>\theta_2>\cdots>\theta_k$, $c_n\downarrow0$, and $\widehat\theta_{j,n}\to\theta_j$ for every $j$.  Then
\[
  \pi_{j,n}(c_n)\to \frac{2(k-j)}{k(k-1)},
  \qquad j=1,\ldots,k.
\]
In particular, when $k>2$, the best treatment has limiting allocation probability $2/k$, not one, and the total limiting allocation probability assigned to inferior treatments is $1-2/k>0$.
\end{proposition}

\begin{proof}
For fixed $j$ and $i$, as $n\to\infty$,
\[
  G\left(\frac{\widehat\theta_{j,n}-\widehat\theta_{i,n}}{c_n}\right)
  \to
  \begin{cases}
  1, & \theta_j>\theta_i,\\
  0, & \theta_j<\theta_i.
  \end{cases}
\]
Treatment $j$ is better than exactly $k-j$ treatments under the strict ordering.  Hence $W_{j,n}(c_n)\to k-j$, and division by $\binom{k}{2}$ gives the result.
\end{proof}

Proposition~\ref{prop:rank-limit-pairwise} is a statement about the conditional probabilities $\pi_{j,n}(c_n)$.  The corresponding statement for the realized \emph{proportions} requires the usual martingale-averaging step and is recorded next.

\begin{corollary}[Proportion limit for direct pairwise annealing]\label{cor:rank-proportions}
Under the hypotheses of Proposition~\ref{prop:rank-limit-pairwise}, on the event $\{\widehat\theta_{j,n}\to\theta_j\ \text{for all }j\}$,
\[
  \frac{N_{j,n}}{n}\to\frac{2(k-j)}{k(k-1)},
  \qquad j=1,\ldots,k,\qquad\text{almost surely.}
\]
\end{corollary}

\begin{proof}
For each $j$, $N_{j,n}=N_{j,km_0}+\sum_{t=km_0+1}^{n}I_{j,t}$ with $\Exp(I_{j,t}\mid\F_{t-1})=\pi_{j,t-1}(c_{t-1})$.  The centered sum $\sum_t\{I_{j,t}-\pi_{j,t-1}(c_{t-1})\}$ is a martingale with increments bounded by $1$, so the bounded-increment martingale SLLN (as in Theorem~\ref{thm:two-proportions}) makes it $o(n)$ almost surely.  On the stated event Proposition~\ref{prop:rank-limit-pairwise} gives $\pi_{j,n}(c_n)\to 2(k-j)/(k(k-1))$, and Ces\`aro averaging of the conditional probabilities, together with $N_{j,km_0}/n\to0$, yields the claim.
\end{proof}

\begin{remark}[On the consistency hypothesis and the worst arm]\label{rem:rank-worst-arm}
Proposition~\ref{prop:rank-limit-pairwise} and Corollary~\ref{cor:rank-proportions} are stated on the consistency event $\{\widehat\theta_{j,n}\to\theta_j\ \forall j\}$, and it is worth recording when this event has probability one.  The non-worst arms are automatically consistent: for $j\le k-1$ the rank limit assigns the strictly positive probability $2(k-j)/(k(k-1))>0$, so $\sum_n\pi_{j,n}(c_n)=\infty$ and, by conditional Borel--Cantelli, arm $j$ is sampled infinitely often, whence $\widehat\theta_{j,n}\to\theta_j$.  Only the worst arm $k$, whose limiting probability is $0$, may be sampled finitely often under fast cooling; on such paths $\widehat\theta_{k,n}$ freezes at a finite value $\widehat\theta_{k,\infty}$.  Even then the rank limit is unaffected, because every comparison involving arm $k$ enters the other arms' scores through the factor $G((\widehat\theta_{i,n}-\widehat\theta_{k,n})/c_n)$, and arm $k$ is abandoned precisely because its score stayed small, i.e.\ $\widehat\theta_{k,\infty}<\theta_i$ for all $i<k$; hence $G((\widehat\theta_{i,n}-\widehat\theta_{k,n})/c_n)\to1$ and $G((\widehat\theta_{k,n}-\widehat\theta_{i,n})/c_n)\to0$ exactly as in the consistent case.  If full consistency of every arm is desired, it can be enforced unconditionally by the forced-exploration device of Section~\ref{sec:forced-exploration}, at the cost of a sublinear number of exploratory allocations.
\end{remark}

\begin{remark}
Thus the direct annealed pairwise-sum design does not recover the Kundu--Jha--Bhandari type concentration on the single best treatment for $k>2$.  It concentrates according to rank scores.  This is a useful correction to the intuitive but false expectation that annealing every pairwise comparison automatically makes the best arm receive asymptotic probability one.
\end{remark}

\subsection{A best-concentrating annealed multi-treatment rule}

To make the best treatment receive probability tending to one, the score must penalize a treatment for losing to even one competitor.  One convenient construction is the product score
\[
  R_{j,n}(c_n)=\prod_{i\neq j}G\left(\frac{\widehat\theta_{j,n}-\widehat\theta_{i,n}}{c_n}\right),
  \qquad
  p_{j,n}=\frac{R_{j,n}(c_n)}{\sum_{\ell=1}^{k}R_{\ell,n}(c_n)}.
\]

\begin{theorem}[Finite inferior allocation for product-score annealing]\label{thm:product-annealed-multi}
Assume $\theta_1>\theta_2\geq\cdots\geq\theta_k$ and let
\[
  \Delta_* = \theta_1-\max_{j\geq2}\theta_j>0.
\]
Suppose that for some $d\in(0,\Delta_*)$ and some almost surely finite $T$,
\[
  \widehat\theta_{1,n}-\widehat\theta_{j,n}\geq d,
  \qquad j=2,\ldots,k,
  \qquad n\geq T.
\]
If
\[
  \sum_{n=1}^{\infty}G(-d/c_n)<\infty,
\]
then the total number of allocations to inferior treatments after time $T$ is finite almost surely and has finite conditional exponential moments.
\end{theorem}

\begin{proof}
For $j\geq2$, the product $R_{j,n}(c_n)$ contains the factor
\[
  G\left(\frac{\widehat\theta_{j,n}-\widehat\theta_{1,n}}{c_n}\right)
  \leq G(-d/c_n).
\]
Also $R_{1,n}(c_n)\geq G(d/c_n)^{k-1}$.  Since $G(d/c_n)\to1$, there is a finite time after which $R_{1,n}(c_n)\geq 2^{-1}$ on the stabilization event.  Hence, for all sufficiently large $n$,
\[
  \sum_{j=2}^{k}p_{j,n}
  \leq
  \frac{\sum_{j=2}^{k}R_{j,n}(c_n)}{R_{1,n}(c_n)}
  \leq 2(k-1)G(-d/c_n).
\]
The result now follows from the same exponential-moment domination argument used in Theorem \ref{thm:annealed-summability}.
\end{proof}

\section{Further Extensions}
\label{sec:further-extensions}

The preceding sections establish the main mechanisms behind fixed and annealed randomized allocation.  Several extensions are now natural.  Clinical-trial data may contain tied best treatments, practically negligible differences, unequal variances, and non-normal responses.  It is also natural to ask how the randomized design behaves against alternatives that are hard to detect.  We take up these questions in turn.

\subsection{Ties and indifference-zone analysis}\label{sec:ext-ties-indifference}

The previous strict-ordering theory should be modified when more than one treatment is truly best or when practically negligible differences are treated as ties.

\subsubsection{Exact ties among the best treatments}

Let
\[
  \theta_1=\cdots=\theta_b=\theta_*>\theta_{b+1}\geq\cdots\geq\theta_k
\]
for some $1<b<k$.  Let $B=\{1,\ldots,b\}$ denote the set of truly best treatments.

\begin{theorem}[Product-score annealing with multiple best treatments]\label{thm:ties-product}
Use the product-score annealed rule of Section \ref{sec:ext-annealed-multi}.  Suppose that for some $d>0$ and some almost surely finite $T$,
\[
  \widehat\theta_{i,n}-\widehat\theta_{j,n}\geq d,
  \qquad i\in B,
  \quad j\notin B,
  \quad n\geq T.
\]
\begin{enumerate}[label=\textnormal{(\alph*)}]
\item \textnormal{(Finite exposure outside $B$.)}  If $\sum_n G(-d/c_n)<\infty$, then the total number of allocations to treatments outside $B$ after $T$ is finite almost surely.
\item \textnormal{(Equal split inside $B$.)}  Suppose, in addition to the consistency $\widehat\theta_{j,n}\to\theta_j$ for all $j$, that the cooling is slow relative to the within-tie estimation error, in the sense that
\begin{equation}\label{eq:slow-cooling-ties}
  \frac{\widehat\theta_{i,n}-\widehat\theta_{i',n}}{c_n}\to 0
  \qquad\text{almost surely, for all } i,i'\in B.
\end{equation}
Then
\[
  p_{j,n}\to \frac1b\quad(j\in B),
  \qquad
  p_{j,n}\to0\quad(j\notin B).
\]
\end{enumerate}
\end{theorem}

\begin{proof}
\emph{(a)}  This is the argument of Theorem \ref{thm:product-annealed-multi}.  For $j\notin B$, the product $R_{j,n}(c_n)$ contains a factor $G((\widehat\theta_{j,n}-\widehat\theta_{i,n})/c_n)$ with $i\in B$, and on the stated event $\widehat\theta_{j,n}-\widehat\theta_{i,n}\le -d$, so this factor is at most $G(-d/c_n)$; meanwhile $R_{i_0,n}(c_n)\ge G(d/c_n)^{k-1}\to1$ for a fixed best arm $i_0$, so eventually $R_{i_0,n}(c_n)\ge1/2$ and $\sum_{j\notin B}p_{j,n}\le 2(k-b)\,G(-d/c_n)$.  Summability and the exponential-moment domination of Theorem \ref{thm:annealed-summability} give a finite outside-$B$ count almost surely.

\emph{(b)}  We compute the limit of each product score $R_{j,n}(c_n)=\prod_{i\neq j}G((\widehat\theta_{j,n}-\widehat\theta_{i,n})/c_n)$.  Split the factors into within-$B$ and across comparisons.

\emph{Step 1 (within-$B$ factors tend to $1/2$).}  This is the step that genuinely needs a hypothesis, because for two tied-best arms $i,i'\in B$ the argument $(\widehat\theta_{i,n}-\widehat\theta_{i',n})/c_n$ is an indeterminate $0/0$ form: the numerator tends to zero (both means converge to $\theta_*$) but so does $c_n$.  Condition \eqref{eq:slow-cooling-ties} states exactly that the numerator is negligible at the current scale, so $(\widehat\theta_{i,n}-\widehat\theta_{i',n})/c_n\to0$ and, by continuity of $G$ at $0$ and symmetry,
\[
  G\!\left(\frac{\widehat\theta_{i,n}-\widehat\theta_{i',n}}{c_n}\right)\to G(0)=\tfrac12,
  \qquad i,i'\in B.
\]

\emph{Step 2 (across factors tend to $0$ or $1$).}  If $j\in B$ and $i\notin B$, then $\widehat\theta_{j,n}-\widehat\theta_{i,n}\ge d>0$ for $n\ge T$, so $(\widehat\theta_{j,n}-\widehat\theta_{i,n})/c_n\to+\infty$ and $G(\cdot)\to1$.  Symmetrically, if $j\notin B$ and $i\in B$, the factor $\to0$.

\emph{Step 3 (scores and normalization).}  For $j\in B$, the product has $b-1$ within-$B$ factors (each $\to1/2$ by Step 1) and $k-b$ across factors against inferior arms (each $\to1$ by Step 2), so
\[
  R_{j,n}(c_n)\to(\tfrac12)^{\,b-1},
  \qquad j\in B,
\]
the same limit for every member of $B$.  For $j\notin B$, the product contains at least one factor against a best arm, which $\to0$ by Step 2, so $R_{j,n}(c_n)\to0$.  Hence $\sum_{\ell=1}^{k}R_{\ell,n}(c_n)\to b\,(\tfrac12)^{b-1}$, and
\[
  p_{j,n}=\frac{R_{j,n}(c_n)}{\sum_{\ell}R_{\ell,n}(c_n)}
  \to\frac{(\tfrac12)^{b-1}}{b\,(\tfrac12)^{b-1}}=\frac1b
  \quad(j\in B),
  \qquad p_{j,n}\to0\quad(j\notin B). \qedhere
\]
\end{proof}

\begin{remark}[A primitive sufficient condition for \eqref{eq:slow-cooling-ties}]\label{rem:slow-cooling-primitive}
By the law of the iterated logarithm, $|\widehat\theta_{i,n}-\theta_*|=O\big(\sqrt{\log\log N_{i,n}/N_{i,n}}\big)$ almost surely for $i\in B$, so a sufficient condition for \eqref{eq:slow-cooling-ties} is
\[
  c_n\sqrt{\frac{N_{i,n}}{\log\log N_{i,n}}}\to\infty
  \qquad\text{a.s. for every } i\in B.
\]
By part (a) all but finitely many allocations fall in $B$, and the within-$B$ symmetry of the rule then makes each best arm receive a linear share $N_{i,n}\asymp n$; for the canonical schedule $c_n=\gamma/\sqrt{\log n}$ this gives $c_n^2 N_{i,n}/\log\log N_{i,n}\asymp n/(\log n\,\log\log n)\to\infty$, so \eqref{eq:slow-cooling-ties} holds.  Crucially, this is compatible with the summability requirement $\sum_n G(-d/c_n)<\infty$ of part (a): for $c_n=\gamma/\sqrt{\log n}$ that requirement reads $\gamma<d/\sqrt2$, while \eqref{eq:slow-cooling-ties} only needs $c_n\sqrt n\to\infty$, which holds for every $\gamma>0$.  Thus the equal-split conclusion and finite outside-$B$ exposure can be achieved simultaneously.
\end{remark}

\subsubsection{Indifference zones}

For a clinically meaningful margin $\eta>0$, define the practically best set
\[
  B_\eta=\{j:\theta_*-\theta_j\leq\eta\},
  \qquad \theta_* = \max_{1\leq j\leq k}\theta_j.
\]
A simple indifference-zone rule is to allocate only among the empirical $\eta$-best treatments
\[
  \widehat B_{\eta,n}
  =\left\{j: \max_{1\leq\ell\leq k}\widehat\theta_{\ell,n}-\widehat\theta_{j,n}\leq\eta\right\},
\]
with optional residual exploration outside this set.

\begin{proposition}[Eventual exclusion outside an indifference zone]\label{prop:indifference-zone}
Assume there is a positive gap beyond the indifference zone, namely
\[
  \min_{j\notin B_\eta}(\theta_* - \theta_j-\eta)>0.
\]
If $\widehat\theta_{j,n}\to\theta_j$ for all $j$, then, almost surely, there is a finite random time $T$ such that
\[
  \widehat B_{\eta,n}\subseteq B_\eta
  \qquad\text{for all }n\geq T.
\]
If the rule assigns probability at most $r_n$ to treatments outside $\widehat B_{\eta,n}$ after time $n$, then the number of allocations outside $B_\eta$ is finite when $\sum_n r_n<\infty$ and is sublinear when $\sum_{n\leq N}r_n=o(N)$ under the martingale normalization condition \(\sum_{N}r_N/R_N^2<\infty\), where \(R_N=\sum_{n\leq N}r_n\).
\end{proposition}

\begin{proof}
Let
\[
  \rho=\min_{j\notin B_\eta}(\theta_* - \theta_j-\eta)>0.
\]
By consistency, for all sufficiently large $n$, every sample mean is within $\rho/3$ of its true mean.  If $j\notin B_\eta$ and $i$ is a best treatment, then
\[
  \widehat\theta_{i,n}-\widehat\theta_{j,n}
  \geq (\theta_* -\rho/3)-(\theta_j+\rho/3)
  =\theta_* -\theta_j-2\rho/3
  >\eta.
\]
Thus $j\notin\widehat B_{\eta,n}$ eventually, proving the set inclusion.  The finite claim follows from conditional Borel--Cantelli when \(\sum_n r_n<\infty\).  The sublinear claim follows from a bounded-increment martingale argument together with the stated normalization condition, exactly as in the proof of Theorem \ref{thm:cooling-classification}.
\end{proof}

\subsection{Unknown and unequal variances}\label{sec:ext-unknown-variances}

The fixed randomized allocation results above do not require equal variances.  Equal-variance normality is convenient for exposition, but the proofs of allocation proportions use only the SLLN for the sample means and the martingale SLLN for treatment indicators.

\begin{theorem}[Mean-based rules with unequal variances]\label{thm:unequal-variances}
Suppose that for each treatment $j$, the observations $X_{j,1},X_{j,2},\ldots$ are independent and identically distributed with mean $\theta_j$ and finite variance $\sigma_j^2\in(0,\infty)$, and suppose the streams are mutually independent.  Then all fixed mean-based randomized allocation proportion results in Theorems \ref{thm:two-proportions} and \ref{thm:multi-proportions} remain valid without assuming $\sigma_0^2=\sigma_1^2$ or known variances.
\end{theorem}

\begin{proof}
The proof of Lemma \ref{lem:both-infinite} uses only the positivity of the conditional allocation probabilities under Assumption \ref{ass:G}.  Once every treatment is sampled infinitely often, the SLLN gives $\widehat\theta_{j,n}\to\theta_j$ along the adaptive sample sizes.  The allocation proportion theorem then follows from the bounded-increment martingale SLLN.  The variances do not enter these steps.
\end{proof}

\begin{proposition}[Consistency of adaptive sample variances]\label{prop:sample-var-consistency}
Under the assumptions of Theorem \ref{thm:unequal-variances}, define
\[
  \widehat\sigma_{j,n}^2
  =\frac{1}{N_{j,n}-1}\sum_{r=1}^{N_{j,n}}(X_{j,r}-\widehat\theta_{j,n})^2
\]
whenever $N_{j,n}\geq2$.  Under any design for which $N_{j,n}\to\infty$ almost surely,
\[
  \widehat\sigma_{j,n}^2\to\sigma_j^2
  \qquad\text{a.s.}
\]
\end{proposition}

\begin{proof}
The sample variance based on the first $m$ observations from treatment $j$ converges almost surely to $\sigma_j^2$ as $m\to\infty$.  Since $N_{j,n}\to\infty$ almost surely, evaluating this convergent sequence at the random indices $N_{j,n}$ preserves the limit.
\end{proof}

\begin{remark}[Studentization changes the effective cooling]
If one replaces a mean contrast by a statistic whose denominator tends to zero, such as
\[
  \frac{\widehat\theta_{0,n}-\widehat\theta_{1,n}}
  {\sqrt{\widehat\sigma_{0,n}^2/N_{0,n}+\widehat\sigma_{1,n}^2/N_{1,n}}},
\]
then the rule is no longer a fixed smooth function of the raw treatment contrast.  It has an implicit sample-size-dependent sharpening.  Such a rule should therefore be analyzed using the annealing principles of Sections \ref{sec:ext-annealed-two} and \ref{sec:ext-cooling-classification}, rather than by applying the fixed-$G$ theorem directly.
\end{remark}

\subsection{General response distributions}\label{sec:ext-general-distributions}

The core allocation-proportion results are not special to normal responses.  Normality becomes useful only when one wants explicit exponential tail bounds and finite moment bounds for stabilization times.

\begin{theorem}[Fixed randomized allocation under an SLLN]\label{thm:general-slln}
For each treatment $j$, suppose the response stream satisfies
\[
  \frac1m\sum_{r=1}^{m}X_{j,r}\to\theta_j
  \qquad\text{a.s.}
\]
as $m\to\infty$.  Under Assumption \ref{ass:G}, Theorems~\ref{thm:two-proportions} and \ref{thm:multi-proportions} remain valid.  Consequently, Corollaries~\ref{cor:two-moments}, \ref{cor:min-count}, and \ref{cor:multi-linear} also remain valid under this SLLN assumption.
\end{theorem}

\begin{proof}
We indicate the two-treatment proof first.  Lemma~\ref{lem:both-infinite} uses only two facts: the conditional allocation probability is adapted to the past, and Assumption~\ref{ass:G} keeps this probability strictly positive whenever the empirical contrast is finite.  Hence the conditional Borel--Cantelli argument still implies $N_{0,n}\to\infty$ and $N_{1,n}\to\infty$ almost surely.  Evaluating the assumed SLLN for the $j$th response stream at the random indices $N_{j,n}$ then gives $\widehat\theta_{j,n}\to\theta_j$ almost surely for $j=0,1$.

With this consistency established, the conditional allocation probability satisfies
\[
  G(\widehat\theta_{0,n}-\widehat\theta_{1,n})\to G(\theta_0-\theta_1)
  \qquad\text{a.s.}
\]
The martingale decomposition used in Theorem~\ref{thm:two-proportions} has bounded increments and is unchanged by the response distribution.  The martingale SLLN and Cesaro averaging therefore prove Theorem~\ref{thm:two-proportions}, and the corresponding corollaries follow exactly as before.  The multi-treatment proof is identical after replacing the two conditional probabilities by the normalized probabilities $\pi_{j,n}$ in \eqref{eq:multi-rule}; positivity gives infinite sampling, the assumed SLLN gives consistency of all adaptive sample means, and martingale averaging gives Theorem~\ref{thm:multi-proportions}.
\end{proof}

For annealed or deterministic bounded-exposure statements, one needs control over the time after which sample means remain correctly ordered.  The following condition is sufficient.

\begin{assumption}[Exponential concentration of sample means]\label{ass:exp-concentration}
For each treatment $j$ and every $u>0$, there exist constants $C_j(u)<\infty$ and $a_j(u)>0$ such that
\[
  \Prob\left(\left|\frac1m\sum_{r=1}^{m}X_{j,r}-\theta_j\right|>u\right)
  \leq C_j(u)e^{-a_j(u)m},
  \qquad m\geq1.
\]
\end{assumption}

This condition holds for many bounded, sub-Gaussian, and regular exponential-family responses.  It includes Bernoulli responses and normal responses with finite variance proxy.

\begin{lemma}[Exponential tail for stabilization time]\label{lem:general-stabilization}
Under Assumption \ref{ass:exp-concentration}, define
\[
  T_u=\inf\left\{M:\left|\frac1m\sum_{r=1}^{m}X_{j,r}-\theta_j\right|\leq u
  \text{ for all }m\geq M\text{ and all }j\right\}.
\]
Then $T_u<\infty$ almost surely and has an exponential tail.  In particular, all moments of $T_u$ are finite.
\end{lemma}

\begin{proof}
By the union bound,
\[
\begin{aligned}
  \Prob(T_u>M)
  &\leq \sum_j\sum_{m=M}^{\infty}
  \Prob\left(\left|\frac1m\sum_{r=1}^{m}X_{j,r}-\theta_j\right|>u\right)\\
  &\leq \sum_j\sum_{m=M}^{\infty}C_j(u)e^{-a_j(u)m}
  \leq C(u)e^{-a(u)M}
\end{aligned}
\]
for suitable constants $C(u)<\infty$ and $a(u)>0$.  The exponential tail implies finiteness almost surely and finiteness of all moments.
\end{proof}

\begin{theorem}[Boundedness mechanisms beyond normality]\label{thm:general-bounded-mechanisms}
Assume a positive treatment gap and Assumption \ref{ass:exp-concentration}.  Any deterministic or annealed rule whose post-stabilization inferior allocation probabilities are summable has a finite post-stabilization inferior allocation count with all moments finite.  If the stabilization time is included in the count, all moments remain finite.
\end{theorem}

\begin{proof}
Lemma \ref{lem:general-stabilization} gives a stabilization time with all moments finite.  After stabilization, the inferior allocation count is dominated by a conditionally Bernoulli sum with summable conditional probabilities.  The exponential-moment argument of Theorem \ref{thm:annealed-summability} gives finite moments for the post-stabilization count.  The sum of two random variables with finite moments of all orders has finite moments of all orders.
\end{proof}

\subsection{Local alternatives: allocation balance and asymptotic power}\label{sec:ext-local-alternatives}

For a \emph{fixed} separation $\Delta=\theta_0-\theta_1>0$ the testing picture is degenerate.  The studentized test of Theorem~\ref{thm:clt-contrast} is consistent, so its power tends to one, while the inferior-allocation fraction tends to the constant $G(-\Delta)>0$ of Corollary~\ref{cor:two-moments}.  Neither limit separates one design from another.  The informative regime is that of \emph{local}, or Pitman, alternatives, where the separation shrinks with the sample size at the rate of the estimation precision.  We treat this regime for a general response distribution and determine, in it, both the limiting inferior-allocation fraction and the limiting power.  The normal model of the earlier sections will reappear as a special case.

We place a mild regularity condition on the responses.

\begin{assumption}[Regular location family]\label{ass:location-family}
For each arm $j=0,1$ the responses are $X_{j,r}=\theta_j+\varepsilon_{j,r}$, where $(\varepsilon_{j,r})_{r\geq1}$ are i.i.d.\ from a density $f_j$ that does not depend on $n$, with mean zero and finite variance $\sigma_j^2\in(0,\infty)$.  The inferior-arm density $f_1$ is absolutely continuous, with score $\dot\ell_1=f_1'/f_1$ and finite Fisher information
\[
  I_1=\int \left(\frac{f_1'(x)}{f_1(x)}\right)^2 f_1(x)\,dx\in(0,\infty),
\]
and $\sqrt{f_1}$ is differentiable in quadratic mean.  When efficiency is discussed, the same regularity is imposed on $f_0$, with score $\dot\ell_0$ and Fisher information $I_0$.
\end{assumption}

The differentiability-in-quadratic-mean condition is the standard regularity for local asymptotic normality \citep[Chapter~7]{vanderVaart1998}.  The normal model is the special case $f_j=N(0,\sigma_j^2)$, for which $\dot\ell_j(x)=-x/\sigma_j^2$ and $I_j=1/\sigma_j^2$.

The local alternative perturbs only the inferior arm,
\begin{equation}\label{eq:local-alt}
  H_{1,n}:\quad \theta_0=\mu,\qquad
  \theta_1^{(n)}=\mu-\frac{h}{\sqrt n},\qquad
  \Delta_n:=\theta_0-\theta_1^{(n)}=\frac{h}{\sqrt n},
\end{equation}
with $h\geq0$ fixed, so $h=0$ recovers the null $H_0:\theta_0=\theta_1=\mu$.  Let $P_{0,n}$ and $P_{1,n}$ denote the laws of the first $n$ assignment--response pairs under $H_0$ and $H_{1,n}$.  We use the fixed rule \eqref{eq:two-rule-0}--\eqref{eq:two-rule-1} with $G$ symmetric, continuous, and strictly increasing at the origin.  The assignment at each stage has conditional probability $G(\pm\widehat\Delta_{t-1})$, a function of the observed past alone.  This last fact will be used twice.

The $n^{-1/2}$ scale is dictated by the precision of the contrast.  The sample-mean contrast is estimated at the rate $n^{-1/2}$ (Theorem~\ref{thm:clt-means}), so a separation is detected with a non-degenerate probability precisely when it is itself of order $n^{-1/2}$.  Proposition~\ref{prop:local-scaling} makes this exact.  We first record that the design balances at the null, and that this holds for any response law.

\begin{lemma}[Balance at the null]\label{lem:balance-null}
Under $H_0$, and for any response distribution with a finite mean, the fixed rule gives $N_{j,n}/n\to\tfrac12$ almost surely for $j=0,1$.
\end{lemma}

\begin{proof}
This is Theorem~\ref{thm:two-proportions} evaluated at $\Delta=0$.  Lemma~\ref{lem:both-infinite} uses only Assumption~\ref{ass:G}, not the gap, so $N_{j,n}\to\infty$.  Lemma~\ref{lem:consistency} then gives $\widehat\theta_{j,n}\to\mu$, whence the conditional probability $G(\widehat\theta_{0,t}-\widehat\theta_{1,t})\to G(0)=\tfrac12$ almost surely.  Ces\`aro averaging of the bounded-increment martingale decomposition of Theorem~\ref{thm:two-proportions} yields $N_{0,n}/n\to\tfrac12$.  The case $j=1$ is symmetric.
\end{proof}

\begin{lemma}[Local asymptotic normality and contiguity]\label{lem:contiguity}
Under Assumption~\ref{ass:location-family} and the sequence \eqref{eq:local-alt}, the log-likelihood ratio $\Lambda_n=\log(dP_{1,n}/dP_{0,n})$ admits the expansion
\[
  \Lambda_n
  =\frac{h}{\sqrt n}\sum_{r=1}^{N_{1,n}}\dot\ell_1(\varepsilon_{1,r})
   -\frac{h^2}{2}\,I_1\,\frac{N_{1,n}}{n}+o_P(1),
\]
and under $H_0$,
\[
  \Lambda_n\dto N\!\left(-\tfrac12 v^2,\;v^2\right),
  \qquad v^2=\tfrac12\,h^2 I_1.
\]
Consequently the sequence $(P_{1,n})$ is contiguous to $(P_{0,n})$.
\end{lemma}

\begin{proof}
The assignment probabilities $G(\pm\widehat\Delta_{t-1})$ are identical functions of the observed data under $H_0$ and $H_{1,n}$, so they cancel in the likelihood ratio.  Arm $0$ has the same law under both hypotheses.  Only the $N_{1,n}$ inferior-arm responses therefore contribute, and under $H_0$, where $X_{1,r}-\mu=\varepsilon_{1,r}$,
\[
  \Lambda_n
  =\sum_{r=1}^{N_{1,n}}
   \log\frac{f_1\!\big(\varepsilon_{1,r}+h/\sqrt n\big)}{f_1(\varepsilon_{1,r})}.
\]
Differentiability in quadratic mean gives, for a deterministic index $m$,
\[
  \sum_{r=1}^{m}\log\frac{f_1(\varepsilon_{1,r}+h/\sqrt n)}{f_1(\varepsilon_{1,r})}
  =\frac{h}{\sqrt n}\sum_{r=1}^{m}\dot\ell_1(\varepsilon_{1,r})
   -\frac{h^2}{2}\,I_1\,\frac{m}{n}+R_{m,n},
\]
with $\sup_{m\leq Kn}|R_{m,n}|\pto0$ for every $K<\infty$ \citep[Chapter~7]{vanderVaart1998}.  The inferior-arm responses form a fixed i.i.d.\ stream, since the rule chooses how many, not which, of its draws are used (Remark~\ref{rem:no-bias}), and $N_{1,n}\leq n$.  Evaluating at $m=N_{1,n}$ gives the stated expansion.  Under $H_0$, $N_{1,n}/n\to\tfrac12$ almost surely, so the quadratic term tends to $-\tfrac14 h^2 I_1=-\tfrac12 v^2$ in probability.  The scores $\dot\ell_1(\varepsilon_{1,r})$ are i.i.d.\ with mean zero and variance $I_1$, and $N_{1,n}/(n/2)\to1$ almost surely, so Anscombe's theorem (Section~\ref{sec:inference}) yields
\[
  \frac{1}{\sqrt n}\sum_{r=1}^{N_{1,n}}\dot\ell_1(\varepsilon_{1,r})
  =\sqrt{\frac{N_{1,n}}{n}}\cdot\frac{1}{\sqrt{N_{1,n}}}\sum_{r=1}^{N_{1,n}}\dot\ell_1(\varepsilon_{1,r})
  \dto \sqrt{\tfrac12}\,N(0,I_1)=N\!\big(0,\tfrac12 I_1\big).
\]
Hence the linear term converges in law to $h\,N(0,\tfrac12 I_1)=N(0,v^2)$, and $\Lambda_n\dto N(-\tfrac12 v^2,v^2)$.  The limiting mean equals $-\tfrac12$ times the limiting variance, so Le Cam's first lemma gives contiguity of $(P_{1,n})$ to $(P_{0,n})$ \citep[Chapter~6]{vanderVaart1998}.  In the normal case $\dot\ell_1(x)=-x/\sigma_1^2$ and $I_1=1/\sigma_1^2$, which recovers $v^2=h^2/(2\sigma_1^2)$.
\end{proof}

\begin{corollary}[Balance under the local alternative]\label{cor:balance-local}
Under the sequence \eqref{eq:local-alt}, $N_{j,n}/n\pto\tfrac12$ for $j=0,1$.  Equivalently, the limiting inferior-allocation fraction is $\tfrac12$: under contiguous local alternatives the fixed randomized design is asymptotically balanced.
\end{corollary}

\begin{proof}
By Lemma~\ref{lem:balance-null}, $N_{1,n}/n-\tfrac12\to0$ almost surely, hence in $P_{0,n}$-probability.  Contiguity (Lemma~\ref{lem:contiguity}) transfers $o_P(1)$ statements from $(P_{0,n})$ to $(P_{1,n})$.  So $N_{1,n}/n-\tfrac12\pto0$ under $H_{1,n}$ as well, and $N_{0,n}/n=1-N_{1,n}/n\pto\tfrac12$.
\end{proof}

We can now give the limiting law of the test statistic and the resulting power.  Recall from Theorem~\ref{thm:clt-contrast} the studentized statistic
\[
  Z_n=\frac{\widehat\theta_{0,n}-\widehat\theta_{1,n}}{s_n},
  \qquad
  s_n^2=\frac{\widehat\sigma_{0,n}^2}{N_{0,n}}+\frac{\widehat\sigma_{1,n}^2}{N_{1,n}}.
\]

\begin{theorem}[Local power of the studentized test]\label{thm:local-power}
Suppose the responses have finite variances $\sigma_0^2,\sigma_1^2$, and let Assumption~\ref{ass:location-family} hold for the inferior arm.  Under the local alternative \eqref{eq:local-alt},
\[
  Z_n\dto N(\delta,1),
  \qquad
  \delta=\frac{h}{\sqrt{2(\sigma_0^2+\sigma_1^2)}}.
\]
Consequently the limiting power of the asymptotically size-$\alpha_0$ tests of $H_0:\theta_0=\theta_1$ is
\[
  \text{one-sided }(Z_n>z_{1-\alpha_0}):\quad
  \beta(h)=\Phi(\delta-z_{1-\alpha_0}),
\]
\[
  \text{two-sided }(|Z_n|>z_{1-\alpha_0/2}):\quad
  \beta(h)=\Phi(\delta-z_{1-\alpha_0/2})+\Phi(-\delta-z_{1-\alpha_0/2}),
\]
where $z_q=\Phi^{-1}(q)$.  At $h=0$ these reduce to $\alpha_0$, which confirms asymptotic size control.
\end{theorem}

\begin{proof}
Write $D_n=\widehat\theta_{0,n}-\widehat\theta_{1,n}$ and decompose
\[
  D_n-\Delta_n=(\widehat\theta_{0,n}-\mu)-(\widehat\theta_{1,n}-\theta_1^{(n)}).
\]
The centred responses are $X_{0,r}-\mu=\varepsilon_{0,r}$ and $X_{1,r}-\theta_1^{(n)}=\varepsilon_{1,r}$, each i.i.d.\ with mean zero and the fixed variance $\sigma_j^2$.  No normality is needed here.  By Corollary~\ref{cor:balance-local}, $N_{j,n}/(n/2)\pto1$, so Anscombe's theorem (Section~\ref{sec:inference}) gives, for each arm,
\[
  \sqrt n\,(\widehat\theta_{j,n}-\theta_j^{(n)})
  =\sqrt{\frac{n}{N_{j,n}}}\cdot\frac{1}{\sqrt{N_{j,n}}}\sum_{r=1}^{N_{j,n}}\varepsilon_{j,r}
  \dto N(0,2\sigma_j^2),
\]
with $\theta_0^{(n)}=\mu$, using $n/N_{j,n}\pto2$ and Slutsky.  The two arms are built from the independent families $(\varepsilon_{0,r})$ and $(\varepsilon_{1,r})$ with asymptotically deterministic normalizers, so the bivariate Anscombe argument of Theorem~\ref{thm:clt-means} makes their limits independent.  Hence
\[
  \sqrt n\,(D_n-\Delta_n)\dto N(0,\tau^2),
  \qquad \tau^2=2(\sigma_0^2+\sigma_1^2).
\]
Since $\sqrt n\,\Delta_n=h$ is constant, $\sqrt n\,D_n\dto N(h,\tau^2)$.  For the denominator, Proposition~\ref{prop:sample-var-consistency} gives $\widehat\sigma_{j,n}^2\pto\sigma_j^2$, the centred residuals having a fixed law, and with $n/N_{j,n}\pto2$,
\[
  n\,s_n^2=\frac{\widehat\sigma_{0,n}^2}{N_{0,n}/n}+\frac{\widehat\sigma_{1,n}^2}{N_{1,n}/n}
  \pto 2\sigma_0^2+2\sigma_1^2=\tau^2 .
\]
Therefore, by Slutsky's theorem,
\[
  Z_n=\frac{\sqrt n\,D_n}{\sqrt{n\,s_n^2}}\dto \frac{N(h,\tau^2)}{\tau}=N\!\left(\frac{h}{\tau},1\right),
\]
which is the claim with $\delta=h/\tau$.  The power expressions follow from $\Prob(N(\delta,1)>z_{1-\alpha_0})=\Phi(\delta-z_{1-\alpha_0})$ and the analogous two-sided computation.  At $h=0$, $\delta=0$ and the right-hand sides equal $\alpha_0$.
\end{proof}

The noncentrality is built from the response variances, not the Fisher informations.  The studentized mean contrast is therefore a distribution-free statistic, and its local power does not know whether the responses are normal.  The price of this robustness is a loss of efficiency away from the normal model, which we quantify in Remark~\ref{rem:efficiency}.

\begin{proposition}[The local scaling is the unique informative one]\label{prop:local-scaling}
Replace \eqref{eq:local-alt} by $\Delta_n=h\,n^{-\rho}$ with $h>0$ and $\rho\geq0$, keeping the fixed rule and finite response variances.  For the size-$\alpha_0$ studentized test the following hold.
\begin{enumerate}[label=\textnormal{(\roman*)}]
\item If $\rho<\tfrac12$, the power tends to $1$, so the test is consistent.
\item If $\rho>\tfrac12$, the power tends to the size $\alpha_0$.
\item If $\rho=\tfrac12$, the power tends to $\beta(h)\in(\alpha_0,1)$ of Theorem~\ref{thm:local-power}.
\end{enumerate}
Thus $\rho=\tfrac12$ is the unique rate at which the limiting power is non-degenerate, which is what singles out the choice \eqref{eq:local-alt}.
\end{proposition}

\begin{proof}
We use the decomposition from the proof of Theorem~\ref{thm:local-power}.  Since $n/N_{j,n}\geq1$, we have $n\,s_n^2\geq\widehat\sigma_{0,n}^2+\widehat\sigma_{1,n}^2\pto\sigma_0^2+\sigma_1^2>0$.

\emph{Case $\rho>\tfrac12$.}  The log-likelihood ratio of Lemma~\ref{lem:contiguity} now has limiting variance $v_n^2=\tfrac12 h^2 n^{1-2\rho}I_1\to0$, so $\Lambda_n\pto0$ and $(P_{1,n})$ is contiguous to $(P_{0,n})$.  Balance and $n\,s_n^2\pto\tau^2$ therefore transfer to $H_{1,n}$, and $\sqrt n\,(D_n-\Delta_n)\dto N(0,\tau^2)$.  As $\sqrt n\,\Delta_n=h\,n^{1/2-\rho}\to0$, we get $Z_n\dto N(0,1)$, and the rejection probability tends to $\alpha_0$.

\emph{Case $\rho=\tfrac12$.}  This is Theorem~\ref{thm:local-power}, with limiting power $\beta(h)\in(\alpha_0,1)$.

\emph{Case $0\leq\rho<\tfrac12$.}  Here $v_n^2=\tfrac12 h^2 n^{1-2\rho}I_1\to\infty$, so the log-likelihood ratio is not tight and contiguity fails.  We argue directly.  Both arms keep a positive limiting share of the allocation.  If $\rho=0$, the gap equals the constant $h$, and Theorem~\ref{thm:two-proportions} gives $N_{0,n}/n\to G(h)$ and $N_{1,n}/n\to G(-h)$, both in $(0,1)$.  If $0<\rho<\tfrac12$, then $\Delta_n\to0$ and the rule is asymptotically symmetric, so $N_{j,n}/n\pto\tfrac12$ by the martingale-averaging argument of Theorem~\ref{thm:two-proportions}, which uses only $\Delta_n\to0$ and not contiguity.  In either case $\liminf_n N_{j,n}/n>0$ in probability, so $n/N_{j,n}=O_P(1)$.  The finite variances then give $n\,s_n^2=O_P(1)$ and $\sqrt n\,(D_n-\Delta_n)=O_P(1)$.  Consequently
\[
  Z_n=\frac{\sqrt n\,\Delta_n+O_P(1)}{\sqrt{n\,s_n^2}}
     =\frac{h\,n^{1/2-\rho}+O_P(1)}{O_P(1)}\pto+\infty,
\]
because the numerator diverges while the denominator is bounded in probability and bounded below by $\sqrt{\sigma_0^2+\sigma_1^2}>0$.  The rejection probability therefore tends to $1$.
\end{proof}

\begin{remark}[Randomization is locally free]\label{rem:local-free}
Corollary~\ref{cor:balance-local} shows that, under contiguous local alternatives, the inferior-allocation fraction tends to $\tfrac12$.  To first order the fixed randomized design therefore behaves like the balanced, equal-allocation design.  The limiting variance $\tau^2=2(\sigma_0^2+\sigma_1^2)$ of Theorem~\ref{thm:local-power} is the balanced two-sample variance constant, and $\delta$ is the balanced-design noncentrality.  The fixed smooth rule thus incurs no first-order loss of local power relative to equal allocation, and it buys allocation concealment at no such cost.  This is the local-power face of the persistent-exploration theme of the paper.  The very exposure $G(-\Delta)$ that is ethically costly for a fixed gap becomes, near the null where the arms are hard to tell apart, a balanced allocation that matches equal allocation in power.
\end{remark}

\begin{remark}[Le Cam's third lemma and the role of the design]\label{rem:third-lemma}
The noncentrality $\delta$ also follows from Le Cam's third lemma.  Under $H_0$ the pair $(Z_n,\Lambda_n)$ is asymptotically bivariate normal with $\mathrm{Cov}(Z_n,\Lambda_n)\to\delta$, so the alternative law of $Z_n$ is the null law shifted by $\delta$ \citep[Chapter~6]{vanderVaart1998}.  The covariance equals $\delta$ for any regular location family, because of the score identity $\Exp\{\varepsilon_{1,r}\,\dot\ell_1(\varepsilon_{1,r})\}=-1$.  The computation is transparent for a deeper reason.  As noted in Lemma~\ref{lem:contiguity}, the adaptive assignment probabilities cancel in $\Lambda_n$.  The rule depends on the unknown parameter only through the data, so it is non-informative for the local parameter, and likelihood inference proceeds as if the realized sample sizes were fixed.  This ignorability underlies post-adaptive estimation theory \citep{MelfiPage2000,HuRosenberger2006}.  A bounded-exposure design, by contrast, starves the inferior arm (Remark~\ref{rem:finite-no-clt}).  It cannot stay balanced under local alternatives, and so it forfeits this local power.
\end{remark}

\begin{remark}[Efficiency away from the normal model]\label{rem:efficiency}
The studentized mean contrast uses the sample means, which are efficient only in the normal location model.  Suppose Assumption~\ref{ass:location-family} holds for both arms.  A procedure built on the efficient (score) estimators replaces each response variance $\sigma_j^2$ by the inverse Fisher information $1/I_j\leq\sigma_j^2$, and attains the noncentrality
\[
  \delta_{\mathrm{eff}}=\frac{h}{\sqrt{2(1/I_0+1/I_1)}}.
\]
The asymptotic relative efficiency of the mean-based test is therefore
\[
  \mathrm{ARE}=\frac{\delta^2}{\delta_{\mathrm{eff}}^2}
  =\frac{1/I_0+1/I_1}{\sigma_0^2+\sigma_1^2}\in(0,1],
\]
with equality if and only if both location families are normal, by the Cram\'er--Rao bound $\sigma_j^2 I_j\geq1$.  The loss of the mean-based test away from normality is thus exactly the classical loss of the sample mean relative to the efficient estimator.  The earlier normal results are recovered by setting $I_j=1/\sigma_j^2$, which gives $v^2=h^2/(2\sigma_1^2)$ and $\mathrm{ARE}=1$.
\end{remark}

\section{Inference after adaptive allocation}
\label{sec:inference}

The previous sections describe \emph{how many} patients each treatment receives.  We now ask what can be \emph{inferred} once the trial stops.  Because the allocation counts $N_{j,n}$ are random and data-dependent, the adaptive sample means are evaluated at random sample sizes, and the relevant limit theory is that of randomly indexed sums.  The qualitative message is that the three exposure regimes of the previous sections---linear, sublinear, and finite inferior exposure---translate into three qualitatively different inferential situations for the inferior treatment, while inference about the superior treatment is always available.

Throughout this section the response streams are i.i.d.\ within each arm with mean $\theta_j$ and finite variance $\sigma_j^2\in(0,\infty)$, and the arms are mutually independent.  We use repeatedly the following randomly indexed central limit theorem.

\begin{proposition*}[Anscombe's random-index CLT]
Let $Y_1,Y_2,\ldots$ be i.i.d.\ with mean $0$ and variance $\sigma^2\in(0,\infty)$, and $S_m=\sum_{r=1}^{m}Y_r$.  Let $(\nu_n)$ be positive integer-valued random variables and $(k_n)$ a deterministic sequence with $k_n\to\infty$ and $\nu_n/k_n\pto1$.  Then
\[
  \frac{S_{\nu_n}}{\sqrt{\nu_n}}\dto N(0,\sigma^2).
\]
The uniform-continuity-in-probability condition required by Anscombe's theorem holds automatically for i.i.d.\ partial sums; see \citet{Anscombe1952} and \citet[Chapter~1]{Gut2009}.  The bivariate form, applied to two independent streams with $\nu_n^{(0)}/k_n^{(0)}\pto1$ and $\nu_n^{(1)}/k_n^{(1)}\pto1$, yields joint convergence to a vector of \emph{independent} normals, because the two coordinate sums are built from independent families and the normalizers are asymptotically deterministic.
\end{proposition*}

\subsection{The fixed-rule (linear) regime}

Under the fixed two-treatment rule, Theorem~\ref{thm:two-proportions} gives $N_{0,n}/n\to\pi_0$ and $N_{1,n}/n\to\pi_1$ almost surely, where
\[
  \pi_0=G(\Delta),\qquad \pi_1=G(-\Delta),\qquad \Delta=\theta_0-\theta_1,\qquad \pi_0+\pi_1=1,
\]
and both limits are strictly positive under Assumption~\ref{ass:G}.  This positivity is exactly what makes large-sample inference possible for both arms.

\begin{theorem}[Asymptotic normality of the adaptive sample means]\label{thm:clt-means}
Under the fixed rule and Assumption~\ref{ass:G}, with finite variances, for each $j\in\{0,1\}$
\[
  \sqrt{N_{j,n}}\,(\widehat\theta_{j,n}-\theta_j)\dto N(0,\sigma_j^2),
  \qquad
  \sqrt{n}\,(\widehat\theta_{j,n}-\theta_j)\dto N\!\left(0,\frac{\sigma_j^2}{\pi_j}\right),
\]
and the two superior/inferior estimators are asymptotically independent:
\[
  \sqrt{n}
  \begin{pmatrix}\widehat\theta_{0,n}-\theta_0\\[2pt]\widehat\theta_{1,n}-\theta_1\end{pmatrix}
  \dto
  N\!\left(\begin{pmatrix}0\\0\end{pmatrix},
  \begin{pmatrix}\sigma_0^2/\pi_0 & 0\\ 0 & \sigma_1^2/\pi_1\end{pmatrix}\right).
\]
\end{theorem}

\begin{proof}
Fix $j$ and set $Y_r=X_{j,r}-\theta_j$, $S_m=\sum_{r=1}^m Y_r$, so that $\widehat\theta_{j,n}-\theta_j=S_{N_{j,n}}/N_{j,n}$.  By Theorem~\ref{thm:two-proportions}, $N_{j,n}\to\infty$ and $N_{j,n}/(\pi_j n)\to1$ almost surely, hence in probability; apply the Anscombe proposition with $\nu_n=N_{j,n}$ and $k_n=\pi_j n$ to get $S_{N_{j,n}}/\sqrt{N_{j,n}}\dto N(0,\sigma_j^2)$, which is the first display.  For the second, write
\[
  \sqrt n(\widehat\theta_{j,n}-\theta_j)
  =\sqrt{\frac{n}{N_{j,n}}}\cdot\frac{S_{N_{j,n}}}{\sqrt{N_{j,n}}},
\]
and since $n/N_{j,n}\to1/\pi_j$ almost surely, Slutsky's theorem gives the variance $\sigma_j^2/\pi_j$.  The joint statement and the asymptotic independence follow from the bivariate Anscombe proposition: the streams $\{X_{0,r}\}$ and $\{X_{1,r}\}$ are independent families, and the normalizers $N_{0,n}/(\pi_0 n)\to1$, $N_{1,n}/(\pi_1 n)\to1$ are asymptotically deterministic, so the limit law factorizes.
\end{proof}

\begin{theorem}[Studentized contrast, confidence interval, and test]\label{thm:clt-contrast}
Maintain the assumptions of Theorem~\ref{thm:clt-means} and let
\[
  s_n^2=\frac{\widehat\sigma_{0,n}^2}{N_{0,n}}+\frac{\widehat\sigma_{1,n}^2}{N_{1,n}},
\]
with $\widehat\sigma_{j,n}^2$ the adaptive sample variances of Proposition~\ref{prop:sample-var-consistency}.  Then
\[
  Z_n:=\frac{(\widehat\theta_{0,n}-\widehat\theta_{1,n})-\Delta}{s_n}\dto N(0,1).
\]
Consequently:
\begin{enumerate}[label=\textnormal{(\roman*)}]
\item the interval $(\widehat\theta_{0,n}-\widehat\theta_{1,n})\pm z_{1-\alpha/2}\,s_n$ is an asymptotically valid $(1-\alpha)$ confidence interval for $\Delta=\theta_0-\theta_1$;
\item the test that rejects $H_0:\theta_0=\theta_1$ when $|(\widehat\theta_{0,n}-\widehat\theta_{1,n})/s_n|>z_{1-\alpha/2}$ has asymptotic level $\alpha$.
\end{enumerate}
\end{theorem}

\begin{proof}
Write the numerator as $(\widehat\theta_{0,n}-\theta_0)-(\widehat\theta_{1,n}-\theta_1)$.  By Theorem~\ref{thm:clt-means} and asymptotic independence,
\[
  \frac{\widehat\theta_{0,n}-\theta_0}{\sqrt{\sigma_0^2/N_{0,n}}}\quad\text{and}\quad
  \frac{\widehat\theta_{1,n}-\theta_1}{\sqrt{\sigma_1^2/N_{1,n}}}
\]
are asymptotically independent standard normals, so the standardized contrast $\{(\widehat\theta_{0,n}-\widehat\theta_{1,n})-\Delta\}/\sqrt{\sigma_0^2/N_{0,n}+\sigma_1^2/N_{1,n}}\dto N(0,1)$.  By Proposition~\ref{prop:sample-var-consistency}, $\widehat\sigma_{j,n}^2\to\sigma_j^2$ almost surely, so $s_n^2/(\sigma_0^2/N_{0,n}+\sigma_1^2/N_{1,n})\to1$, and Slutsky's theorem gives $Z_n\dto N(0,1)$.  Claim (i) is immediate.  For (ii), note that under $H_0$ the two means are equal, so $\Delta=0$ and, by symmetry of $G$, $\pi_0=\pi_1=G(0)=1/2$; both counts still grow linearly ($N_{j,n}/n\to1/2$), so Theorem~\ref{thm:clt-means} applies with $\Delta=0$ and $Z_n=(\widehat\theta_{0,n}-\widehat\theta_{1,n})/s_n\dto N(0,1)$, yielding asymptotic level $\alpha$.
\end{proof}

Theorem~\ref{thm:clt-contrast} is the precise sense in which the fixed randomized rule \emph{buys} valid inference with its persistent exploration: because both $\pi_0$ and $\pi_1$ are bounded away from zero, the contrast is estimated at the parametric $n^{-1/2}$ rate and the usual Wald inference is asymptotically exact.

\subsection{The sublinear and finite regimes}

When the rule is annealed so as to suppress inferior exposure, the inferior estimator is built from far fewer observations, and the inferential guarantee degrades in a way that mirrors the exposure regime.

\begin{proposition}[Inference in the sublinear regime]\label{prop:clt-sublinear}
Consider the forced-exploration rule \eqref{eq:guarded-rule} of Theorem~\ref{thm:forced-exploration} with $\varepsilon_n\downarrow0$, $\sum_n\varepsilon_n=\infty$, $E_N=\sum_{n\le N}\varepsilon_n=o(N)$, and summable cooling.  Then $N_{1,n}\sim\tfrac12 E_n$ almost surely, and
\[
  \sqrt{N_{1,n}}\,(\widehat\theta_{1,n}-\theta_1)\dto N(0,\sigma_1^2),
  \qquad
  \sqrt{\tfrac12 E_n}\,(\widehat\theta_{1,n}-\theta_1)\dto N(0,\sigma_1^2).
\]
The superior estimator continues to satisfy $\sqrt n(\widehat\theta_{0,n}-\theta_0)\dto N(0,\sigma_0^2)$.  Hence inference about $\theta_1$ remains valid but at the slower effective sample size $N_{1,n}\asymp E_n$: a confidence interval for $\theta_1$ has width of order $E_n^{-1/2}$, against $n^{-1/2}$ for $\theta_0$.
\end{proposition}

\begin{proof}
The equivalence $N_{1,n}\sim\tfrac12 E_n$ a.s.\ is Theorem~\ref{thm:forced-exploration}(ii), and consistency $\widehat\theta_{1,n}\to\theta_1$ holds by Theorem~\ref{thm:forced-exploration}(i).  Since $N_{1,n}/(\tfrac12 E_n)\to1$ almost surely, the Anscombe proposition with $\nu_n=N_{1,n}$ and deterministic $k_n=\tfrac12 E_n\to\infty$ gives $S_{N_{1,n}}/\sqrt{N_{1,n}}\dto N(0,\sigma_1^2)$, i.e.\ the first display; the second follows from $\sqrt{N_{1,n}/(\tfrac12E_n)}\to1$ and Slutsky.  The superior arm has $N_{0,n}\sim n$ by the same theorem, and the Anscombe argument of Theorem~\ref{thm:clt-means} applies with $k_n=n$.  The stated interval widths are read off from the normalizing rates.
\end{proof}

\begin{remark}[The finite regime: the statistical price of bounded exposure]\label{rem:finite-no-clt}
In the finite-exposure regime---the summable-cooling annealed rule of Theorem~\ref{thm:annealed-summability}, or any rule with $\sum_n\Prob(A_{n+1}=1\mid\F_n)<\infty$ after stabilization---the inferior treatment receives only finitely many observations almost surely, $N_{1,\infty}<\infty$.  Then $\widehat\theta_{1,n}$ is eventually frozen at the average of those finitely many responses and does \emph{not} converge to $\theta_1$; no law of large numbers, and a fortiori no central limit theorem, is available for the inferior mean.  Inference is therefore possible only for the superior treatment (for which $N_{0,n}\sim n$ and $\sqrt n(\widehat\theta_{0,n}-\theta_0)\dto N(0,\sigma_0^2)$).  This is the exact statistical counterpart of the bounded-exposure property.  The very mechanism that protects patients from the inferior treatment also destroys the information needed to estimate its effect.  Quantitatively, securing even a logarithmic number of inferior observations, of order $\log N$, suffices for a slowly shrinking confidence interval for $\theta_1$, and this already forces the design out of the strictly finite regime.  This is the trade-off made precise.
\end{remark}

\begin{remark}[Multi-treatment inference]\label{rem:multi-inference}
The same analysis applies arm-by-arm in the multi-treatment design of Section~\ref{sec:multi}.  Under the fixed normalized rule, Theorem~\ref{thm:multi-proportions} gives $N_{j,n}/n\to\pi_j^*>0$ for every $j$, so $\sqrt{N_{j,n}}(\widehat\theta_{j,n}-\theta_j)\dto N(0,\sigma_j^2)$ for each arm and any fixed contrast $\theta_a-\theta_b$ admits the Wald inference of Theorem~\ref{thm:clt-contrast} with $(\pi_a^*,\pi_b^*)$ in place of $(\pi_0,\pi_1)$; the arms are asymptotically independent by the multivariate Anscombe argument.  Under product-score annealing (Theorem~\ref{thm:product-annealed-multi}) the inferior arms fall into the finite regime and only the best arm enjoys a CLT, exactly as in Remark~\ref{rem:finite-no-clt}.
\end{remark}

\section{Simulation Study}
\label{sec:simulation-study}

To examine the finite-sample behaviour of the proposed randomized adaptive allocation procedures, we conduct a simulation study under Normal response distributions. The principal objective is to investigate the extent to which the empirical allocation proportions approach their theoretical limits as the total sample size increases. We also examine how the randomization parameter influences the allocation of experimental units across treatments and how the variability of the resulting allocation proportions changes with the sample size. The study considers both the two-treatment rule developed in Section \ref{sec:two} and the normalized multi-treatment rule introduced in Section \ref{sec:multi}.

Throughout the simulations, the response distributions are taken to be independent Normal distributions with treatment-specific means and unit variances. The initial sample size is fixed at one observation per treatment, and the total sample size is varied over $N\in\{25,50,75,100,200,400\}$. The randomization parameter is chosen from $c\in\{1,3,5\}$. For each combination of the treatment parameters, randomization parameter, and total sample size, we perform $10,000$ independent replications. In each replication, the initial observations are generated independently from the specified treatment distributions. Subsequent assignments are determined by the relevant adaptive allocation rule, with each newly assigned experimental unit contributing an independent response from the distribution corresponding to its allocated treatment.

For each configuration, we report the empirical mean and standard deviation of the relevant allocation proportion across the simulation replications. These quantities are compared with the corresponding theoretical limiting allocation proportions derived in Sections \ref{sec:two} and \ref{sec:multi}. In the tables, each entry is presented in the form $a\,(b)$, where $a$ denotes the empirical mean and $b$ denotes the empirical standard deviation of the allocation proportion.

\subsection{Two-treatment randomized allocation}
\label{subsec:two-treatment-simulation}

We consider independent responses $X_{j,r}\sim N(\theta_j,1)$, $j=0,1$, with $\theta_1=0$ and
\[
\theta_0\in\{0,0.3,0.6,0.9,1.2,1.5\}.
\]
One initial observation is assigned to each treatment. The allocation rule assigns the next experimental unit to treatment $0$ with conditional probability given by the probit rule introduced in Section \ref{sec:preliminaries}. Here, $c > 0$ is the randomization parameter, and the complementary probability determines the allocation to treatment $1$.

Let $N_{1,N}$ denote the number of allocations to treatment $1$ by the end of a trial of total size $N$. By Theorem \ref{thm:two-proportions}, when $\theta_0>\theta_1$,
\[
\frac{N_{1,N}}{N}
\xrightarrow{\mathrm{a.s.}}
\Phi\!\left(-\frac{\theta_0-\theta_1}{c}\right).
\]

\autoref{tab:two-treatment-simulation} presents the empirical means and standard deviations of $N_{1,N}/N$ for the specified combinations of $\theta_0$, $c$, and $N$, together with the corresponding theoretical limiting proportions.


The simulation results exhibit a clear relationship between the treatment-mean difference and the allocation of experimental units. For a fixed value of $c$, increasing $\theta_0-\theta_1$ reduces the theoretical limiting proportion allocated to treatment $1$. This behaviour is consistent with the allocation mechanism, which increasingly favours treatment $0$ as the difference between the treatment means becomes more pronounced.


The effect of the randomization parameter is also apparent. For a fixed positive treatment-mean difference, increasing $c$ moves the theoretical limiting allocation proportion of the inferior treatment closer to $1/2$. Conversely, a smaller value of $c$ produces a more pronounced allocation preference for the treatment with the larger mean. For example, when $\theta_0-\theta_1=0.9$, the theoretical limiting proportion allocated to treatment $1$ is approximately $0.1841$, $0.3821$, and $0.4286$ for $c=1,3,$ and $5$, respectively. These values illustrate the extent to which the randomization scale regulates the balance between preferential allocation and continued exploration.

The empirical standard deviations generally decrease as the total sample size increases. This pattern indicates that the realized allocation proportions become less variable across independent replications as the trial progresses. Moreover, for a fixed treatment configuration, the empirical means tend to approach their corresponding theoretical limits as $N$ increases, although small finite-sample fluctuations remain. Taken together, these findings provide a numerical illustration of the allocation-proportion convergence established in Theorem \ref{thm:two-proportions}.

An important implication is that a smaller randomization scale can substantially reduce the long-run proportion assigned to the inferior treatment, but it does not eliminate that proportion under fixed non-degenerate randomization. In accordance with Corollary \ref{cor:two-moments}, the limiting inferior-treatment proportion remains strictly positive whenever the treatment means are distinct and $c>0$. The simulation results therefore illustrate that, under fixed probit randomization, the inferior treatment retains a strictly positive limiting allocation proportion, even when the randomization parameter is small.

\begin{table}[!htbp]
\centering
\caption{
Empirical means of $N_{1,N}/N$ under the fixed probit randomized allocation rule.
The corresponding standard deviations are reported in parentheses.
The last column gives the theoretical limiting allocation proportion
$\Phi\!\left(-(\theta_0-\theta_1)/c\right)$.
}
\label{tab:two-treatment-simulation}

\small
\setlength{\tabcolsep}{4pt}
\renewcommand{\arraystretch}{1.15}

\resizebox{\textwidth}{!}{
\begin{tabular}{cccccccc}

\toprule

\multirow{2}{*}{$\theta_0$}
&
\multicolumn{6}{c}{$E(N_{1,N}/N)\;(\mathrm{sd}(N_{1,N}/N))$}
&
\multirow{2}{*}{$\Phi\!\left(-\dfrac{\theta_0-\theta_1}{c}\right)$}
\\

\cmidrule(lr){2-7}

&
$N=25$
&
$N=50$
&
$N=75$
&
$N=100$
&
$N=200$
&
$N=400$
&
\\

\midrule

\multicolumn{8}{c}{\textbf{$c=1$}}
\\

\midrule

0.0 &
0.5055 (0.2382) &
0.5095 (0.2029) &
0.5015 (0.1832) &
0.4897 (0.1632) &
0.5001 (0.1122) &
0.5039 (0.0776) &
0.5000
\\

0.3 &
0.4058 (0.2373) &
0.3855 (0.1954) &
0.3786 (0.1771) &
0.3748 (0.1542) &
0.3805 (0.1156) &
0.3787 (0.0764) &
0.3821
\\

0.6 &
0.2998 (0.2072) &
0.2891 (0.1754) &
0.2708 (0.1585) &
0.2700 (0.1385) &
0.2648 (0.1119) &
0.2660 (0.0813) &
0.2743
\\

0.9 &
0.2402 (0.1839) &
0.2019 (0.1454) &
0.1837 (0.1297) &
0.1785 (0.1099) &
0.1743 (0.0867) &
0.1750 (0.0633) &
0.1841
\\

1.2 &
0.1834 (0.1422) &
0.1398 (0.1191) &
0.1246 (0.0953) &
0.1218 (0.0859) &
0.1100 (0.0668) &
0.1042 (0.0512) &
0.1151
\\

1.5 &
0.1326 (0.1122) &
0.0981 (0.0806) &
0.0837 (0.0726) &
0.0754 (0.0613) &
0.0652 (0.0469) &
0.0618 (0.0372) &
0.0668
\\

\midrule

\multicolumn{8}{c}{\textbf{$c=3$}}
\\

\midrule

0.0 &
0.4975 (0.1189) &
0.5034 (0.0866) &
0.5049 (0.0685) &
0.4999 (0.0644) &
0.4990 (0.0433) &
0.4982 (0.0312) &
0.5000
\\

0.3 &
0.4583 (0.1183) &
0.4631 (0.0854) &
0.4623 (0.0747) &
0.4618 (0.0629) &
0.4596 (0.0447) &
0.4606 (0.0306) &
0.4602
\\

0.6 &
0.4221 (0.1178) &
0.4247 (0.0875) &
0.4237 (0.0717) &
0.4229 (0.0624) &
0.4185 (0.0449) &
0.4200 (0.0309) &
0.4207
\\

0.9 &
0.3812 (0.1148) &
0.3838 (0.0859) &
0.3817 (0.0705) &
0.3855 (0.0607) &
0.3820 (0.0422) &
0.3837 (0.0316) &
0.3821
\\

1.2 &
0.3604 (0.1141) &
0.3496 (0.0838) &
0.3493 (0.0720) &
0.3471 (0.0580) &
0.3453 (0.0432) &
0.3449 (0.0300) &
0.3446
\\

1.5 &
0.3242 (0.1109) &
0.3174 (0.0801) &
0.3136 (0.0665) &
0.3075 (0.0609) &
0.3087 (0.0415) &
0.3097 (0.0287) &
0.3085
\\

\midrule

\multicolumn{8}{c}{\textbf{$c=5$}}
\\

\midrule

0.0 &
0.4983 (0.1076) &
0.4972 (0.0781) &
0.4987 (0.0627) &
0.4971 (0.0565) &
0.5008 (0.0385) &
0.5021 (0.0281) &
0.5000
\\

0.3 &
0.4755 (0.1053) &
0.4771 (0.0749) &
0.4770 (0.0623) &
0.4776 (0.0549) &
0.4762 (0.0396) &
0.4767 (0.0269) &
0.4761
\\

0.6 &
0.4621 (0.1023) &
0.4544 (0.0731) &
0.4571 (0.0624) &
0.4550 (0.0550) &
0.4524 (0.0402) &
0.4531 (0.0278) &
0.4522
\\

0.9 &
0.4361 (0.1052) &
0.4314 (0.0741) &
0.4318 (0.0602) &
0.4326 (0.0517) &
0.4309 (0.0385) &
0.4278 (0.0273) &
0.4286
\\

1.2 &
0.4160 (0.1013) &
0.4100 (0.0730) &
0.4096 (0.0649) &
0.4061 (0.0528) &
0.4089 (0.0381) &
0.4036 (0.0278) &
0.4052
\\

1.5 &
0.3886 (0.0970) &
0.3862 (0.0749) &
0.3834 (0.0600) &
0.3852 (0.0528) &
0.3829 (0.0375) &
0.3818 (0.0266) &
0.3821
\\

\bottomrule

\end{tabular}
}

\end{table}

\subsection{Multi-treatment randomized allocation}
\label{subsec:multi-treatment-simulation}

We consider three treatments with independent responses
$X_{j,r}\sim N(\theta_j,1)$, $j=1,2,3$, and
$\theta_1>\theta_2>\theta_3$. The parameter configurations are
$(\theta_1,\theta_2,\theta_3)\in\{(0.9,0.2,0),(2,1.2,0.5)\}$.
One initial observation is assigned to each treatment. For three treatments ($k=3$), the allocation rule assigns the next experimental unit according to the normalized pairwise probit rule defined in Section \ref{sec:multi}, equation \ref{eq:multi-rule}. With $c>0$ controlling the degree of randomization, the normalization ensures that the allocation probabilities sum to one across the three treatments.

Let $N_{j,N}$ denote the number of allocations to treatment $j$ by the end of a trial of total size $N$. The theoretical limiting allocation proportion, as established in Theorem \ref{thm:multi-proportions}, is
\[
\frac{N_{j,N}}{N}
\xrightarrow{\mathrm{a.s.}}\pi_j^*,
\qquad
\pi_j^*=\frac{1}{3}\sum_{\substack{i=1\\i\ne j}}^3
\Phi\!\left(\frac{\theta_j-\theta_i}{c}\right),
\qquad j=1,2,3.
\]
Since the first treatment has the largest population mean and the second treatment has the second-largest population mean, we focus on $N_{2,N}/N$, the allocation proportion of the second-best treatment. This choice permits a direct assessment of the theoretical allocation behaviour for a treatment that is inferior to the best-performing treatment but superior to the remaining treatment. \autoref{tab:multi-treatment-simulation} reports the empirical means and standard deviations of this proportion, alongside the corresponding theoretical limits.

For the first configuration, $(\theta_1,\theta_2,\theta_3)=(0.9,0.2,0)$, the theoretical limiting proportion for treatment $2$ is
\[
\pi_2^*
=
\frac{1}{3}
\left\{
\Phi\!\left(-\frac{0.7}{c}\right)
+
\Phi\!\left(\frac{0.2}{c}\right)
\right\}.
\]
The resulting limits are approximately $0.2737$, $0.3114$, and $0.3201$ for $c=1,3,$ and $5$, respectively. For the second configuration, $(\theta_1,\theta_2,\theta_3)=(2,1.2,0.5)$, the corresponding expression is
\[
\pi_2^*
=
\frac{1}{3}
\left\{
\Phi\!\left(-\frac{0.8}{c}\right)
+
\Phi\!\left(\frac{0.7}{c}\right)
\right\},
\]
yielding approximate limits of $0.3233$, $0.3290$, and $0.3307$ for the same values of $c$.

The empirical results in \autoref{tab:multi-treatment-simulation} are broadly consistent with the theoretical predictions. Across both treatment configurations, the empirical mean allocation proportions remain close to their respective theoretical limits, with the agreement generally becoming more pronounced at larger sample sizes. The empirical standard deviations also decrease as $N$ increases, reflecting diminishing variability in the allocation proportion across replications. These observations provide finite-sample support for the almost-sure convergence result established in Theorem \ref{thm:multi-proportions}.

The influence of the randomization parameter is particularly
informative. As $c$ increases, the normalized allocation
probabilities become less sensitive to differences between the
current sample means, causing the limiting allocation proportions to become more similar.

The results also illustrate an important distinction between treatment ranking and exposure reduction. In particular, the second-best treatment continues to receive a positive limiting fraction of the total sample, in agreement with Proposition \ref{prop:ordering} and Corollary \ref{cor:multi-linear}. Consequently, the normalized multi-treatment rule maintains continued exploration across all treatments, with each inferior treatment receiving a strictly positive limiting allocation proportion, while the best treatment receives the largest limiting allocation proportion.

\begin{table}[!htbp]
\centering
\caption{
Empirical means of $N_{2,N}/N$ under the normalized fixed probit randomized allocation rule for the three-treatment normal model
($\sigma_1,\sigma_2,\sigma_3)=(1,1,1)$.
The corresponding standard deviations are reported in parentheses.
The last column gives the theoretical limiting allocation proportion
$\pi_2^{*}$.
}
\label{tab:multi-treatment-simulation}

\small
\setlength{\tabcolsep}{4pt}
\renewcommand{\arraystretch}{1.15}

\resizebox{\textwidth}{!}{
\begin{tabular}{cccccccc}

\toprule

\multirow{2}{*}{$(\theta_1,\theta_2,\theta_3)$}
&
\multicolumn{6}{c}{$E(N_{2,N}/N)\;(\mathrm{sd}(N_{2,N}/N))$}
&
\multirow{2}{*}{$\pi_2^{*}$}
\\

\cmidrule(lr){2-7}

&
$N=25$
&
$N=50$
&
$N=75$
&
$N=100$
&
$N=200$
&
$N=400$
&
\\

\midrule

\multicolumn{8}{c}{\textbf{$c=1$}}
\\
\midrule

$(0.9,\,0.2,\,0)$
&
0.2775 (0.1434)
&
0.2750 (0.1263)
&
0.2858 (0.1131)
&
0.2733 (0.1021)
&
0.2780 (0.0783)
&
0.2768 (0.0556)
&
0.2737
\\

$(2,\,1.2,\,0.5)$
&
0.3195 (0.1391)
&
0.3221 (0.1118)
&
0.3247 (0.0992)
&
0.3296 (0.0894)
&
0.3269 (0.0691)
&
0.3262 (0.0482)
&
0.3233
\\

\midrule

\multicolumn{8}{c}{\textbf{$c=3$}}
\\
\midrule

$(0.9,\,0.2,\,0)$
&
0.3118 (0.0987)
&
0.3123 (0.0739)
&
0.3140 (0.0599)
&
0.3143 (0.0547)
&
0.3112 (0.0363)
&
0.3133 (0.0279)
&
0.3114
\\

$(2,\,1.2,\,0.5)$
&
0.3290 (0.0992)
&
0.3291 (0.0764)
&
0.3299 (0.0623)
&
0.3288 (0.0514)
&
0.3279 (0.0373)
&
0.3282 (0.0267)
&
0.3290
\\

\midrule

\multicolumn{8}{c}{\textbf{$c=5$}}
\\
\midrule

$(0.9,\,0.2,\,0)$
&
0.3244 (0.0938)
&
0.3171 (0.0676)
&
0.3215 (0.0566)
&
0.3214 (0.0479)
&
0.3196 (0.0357)
&
0.3191 (0.0239)
&
0.3201
\\

$(2,\,1.2,\,0.5)$
&
0.3295 (0.0966)
&
0.3331 (0.0661)
&
0.3312 (0.0553)
&
0.3289 (0.0496)
&
0.3305 (0.0343)
&
0.3301 (0.0242)
&
0.3307
\\

\bottomrule

\end{tabular}
}
\end{table}

\subsection{Summary of simulation findings}
\label{subsec:simulation-summary}

The two-treatment and three-treatment experiments provide a coherent numerical illustration of the theoretical allocation results developed in this paper. In both settings, the empirical allocation proportions exhibit behaviour consistent with their corresponding theoretical limits, while their variability generally diminishes as the total sample size increases. The randomization parameter $c$ controls the extent to which the allocation mechanism responds to differences between treatment means: smaller values induce stronger preferential allocation, whereas larger values yield more balanced allocation proportions.

Nevertheless, the simulations also underscore the fundamental limitation of fixed non-degenerate randomization. Even when the best treatment is clearly separated from the inferior treatments, the theoretical limiting allocation proportion of every inferior treatment remains strictly positive. Consequently, increasing the total sample size improves the stability of the allocation proportions but does not remove persistent inferior-treatment exposure. These numerical findings are consistent with the impossibility results established in Sections \ref{sec:two} and \ref{sec:multi}, illustrating the limitations of fixed randomized allocation in achieving bounded exposure to inferior treatments.

\section{Real Data Analysis}
\label{sec:realdata}

To illustrate the practical applicability of the proposed normalized multi-treatment randomized adaptive allocation procedure, we consider
the pregabalin clinical trial for the treatment of post-herpetic neuralgia (PHN) reported by \citet{SABATOWSKI2004}. The original study was a
randomized, double-blind, placebo-controlled clinical trial designed to investigate the efficacy and safety of pregabalin at two dosage levels, $150$ mg/day and $300$ mg/day, compared with placebo. Pain intensity was measured using an $11$-point Numerical Rating Scale, with lower endpoint pain scores indicating better pain relief (greater treatment efficacy).

Based on the intention-to-treat (ITT) results, the treatment groups are defined as follows:
\begin{itemize}
    \item Treatment 1: Pregabalin $300$ mg/day, with endpoint least-squares mean pain score $4.76$, standard error $0.23$, and sample size $76$.
    \item Treatment 2: Pregabalin $150$ mg/day, with endpoint least-squares mean pain score $5.14$, standard error $0.22$, and sample size $81$.
    \item Treatment 3: Placebo, with endpoint least-squares mean pain score $6.33$, standard error $0.22$, and sample size $81$.
\end{itemize}

Since individual patient-level observations are unavailable, we use these published summary statistics to construct a model-based redesign of the original trial. To conform to the convention that larger population means correspond to more favourable treatments, the negatives of the reported endpoint pain scores are used as the treatment parameters.

To incorporate clinical responses into our proposed framework, the three treatment responses are assumed to follow independent Normal distributions. Consequently, treatment-specific distributions are specified as follows, 
\[
X_{j,r}\sim N(\theta_j,\sigma_j^2),\qquad j=1,2,3,
\]
\[ \text{with \qquad }
(\theta_1,\theta_2,\theta_3)=(-4.76,-5.14,-6.33),
\qquad \text{and \qquad}
(\sigma_1,\sigma_2,\sigma_3)=(2.00,1.98,1.98),
\]
where treatment $j$ corresponds to the $j$th treatment group described above. The working standard deviations are approximated by $\sigma_j=\mathrm{SE}_j\sqrt{n_j}$ and rounded to two decimal places. The ordering $\theta_1>\theta_2>\theta_3$ identifies pregabalin $300$ mg/day as the best-performing treatment, followed by pregabalin $150$ mg/day and placebo.

Using these true parameter values, we apply the proposed allocation rule described in Section \ref{sec:multi}, assigning one initial observation to each treatment. For each combination of total sample size $N\in\{25,50,75,100,200,400\}$ and randomization scale $c\in\{1,3,5\}$, we perform $10,000$ independent replications.


The performance of the procedure is assessed through the allocation proportion $N_{2,N}/N$, where $N_{2,N}$ denotes the final allocation count for treatment $2$, the second-best treatment. \autoref{Table3} reports its empirical mean and standard deviation across replications. Each entry of the form $a\,(b)$ represents the empirical mean $a=E(N_{2,N}/N)$, followed in parentheses by the empirical standard deviation $\mathrm{sd}(N_{2,N}/N)$. The last column gives the theoretical limiting allocation proportion $\pi_2^*$ derived in the preceding theoretical analysis. For the present parameter configuration, the corresponding theoretical limiting allocation proportion, established in Theorem \ref{thm:multi-proportions}, is
\[
\pi_2^*
=\frac{1}{3}\left[
\Phi\left(\frac{\theta_2-\theta_1}{c}\right)
+
\Phi\left(\frac{\theta_2-\theta_3}{c}\right)
\right]
=\frac{1}{3}\left[
\Phi\left(\frac{-0.38}{c}\right)
+
\Phi\left(\frac{1.19}{c}\right)
\right].
\]
The corresponding theoretical limiting proportions are approximately $0.4116$, $0.3679$, and $0.3546$ for $c=1$, $3$, and $5$, respectively.


As shown in \autoref{Table3}, for each value of $c$, the empirical mean allocation proportions tend towards their theoretical limiting values as the total sample size increases. For $c=1$, the empirical mean rises from $0.3752$ at $N=25$ to $0.4043$ at $N=400$, approaching the limit $0.4116$. For $c=3$, the empirical means at $N=200$ and $N=400$ are both $0.3690$, close to the theoretical value $0.3679$. Similarly, for $c=5$, the empirical mean at $N=400$ is $0.3549$, close to the theoretical limit $0.3546$. These findings are consistent with the convergence of the allocation proportions established in Theorem \ref{thm:multi-proportions} for the considered parameter configuration. Furthermore, the empirical standard deviations also decrease with increasing total sample size for all three values of $c$, indicating reduced variability in the allocation proportion across replications.




Overall, this pregabalin-based analysis provides a clinically motivated illustration of the proposed three-treatment allocation procedure. The close agreement between the empirical allocation proportions and their theoretical limits supports the relevance of the asymptotic allocation results in finite samples.



\begin{table}[!htbp]
\centering
\caption{Results for Pregabalin Drug Trial for Postherpetic Neuralgia (PHN).}
\label{Table3}

\small
\setlength{\tabcolsep}{4pt}
\renewcommand{\arraystretch}{1.15}

\resizebox{\textwidth}{!}{
\begin{tabular}{c|cccccc|c}

\toprule

\multirow{2}{*}{$c$}
&
\multicolumn{6}{c|}{$E(N_{2,N}/N)\;(\mathrm{sd}(N_{2,N}/N))$}
&
\multirow{2}{*}{$\pi_2^{*}$}
\\

\cmidrule(lr){2-7}

&
$N=25$
&
$N=50$
&
$N=75$
&
$N=100$
&
$N=200$
&
$N=400$
&
\\

\midrule

$1$
&
0.3752 (0.1808)
&
0.3825 (0.1675)
&
0.3892 (0.1537)
&
0.3906 (0.1497)
&
0.3976 (0.1250)
&
0.4043 (0.1079)
&
0.4116
\\

\midrule

$3$
&
0.3623 (0.1251)
&
0.3639 (0.0991)
&
0.3669 (0.0837)
&
0.3674 (0.0728)
&
0.3690 (0.0521)
&
0.3690 (0.0362)
&
0.3679
\\

\midrule

$5$
&
0.3506 (0.1048)
&
0.3529 (0.0792)
&
0.3538 (0.0652)
&
0.3539 (0.0562)
&
0.3535 (0.0399)
&
0.3549 (0.0284)
&
0.3546
\\

\bottomrule

\end{tabular}
}
\end{table}

\section{Concluding Remarks}
\label{sec:conclusion}

This paper studied a randomized counterpart of the bounded-exposure design of \citet{KunduJhaBhandari2025}.  For two treatments, the next patient goes to treatment $0$ with probability $G(\widehat\theta_{0,n}-\widehat\theta_{1,n})$ and to treatment $1$ with the complementary probability.  The limiting behaviour is wholly different from that of the deterministic larger-mean rule.  Under any fixed non-degenerate symmetric $G$, both treatments are sampled infinitely often, both sample means are strongly consistent, and the inferior treatment draws the positive limiting fraction $G(\theta_1-\theta_0)$ whenever $\theta_0>\theta_1$.  The inferior count therefore grows linearly, and its moments are not uniformly bounded.

The multi-treatment setting tells the same story, once the pairwise probabilities are correctly normalized.  Every treatment draws a positive limiting fraction, and both the total inferior count and the second-largest count grow linearly.

The two results are thus not rival proofs of one phenomenon.  They describe different designs with different asymptotic ethics.  A deterministic, or progressively sharpened, rule can bound the long-run exposure to an inferior treatment.  A fixed smooth randomized rule keeps randomization and concealment, and pays for them with continued positive exposure.

Annealed randomization interpolates between the two.  By cooling the scale with the accumulated information, it can drive the inferior exposure from linear to sublinear to finite, and the regime is settled by the summability of the tail probabilities.  Annealing alone rests on identifying the correct ordering.  A small forced-exploration floor removes this caveat (Section~\ref{sec:forced-exploration}), since it makes correct identification almost sure while keeping the exposure sublinear.  The inference of Section~\ref{sec:inference} then shows that the three regimes are far from interchangeable.  Linear exposure gives the usual root-$n$ Wald inference for the inferior treatment, sublinear exposure gives valid but slower inference at the exploration rate, and finite exposure allows no central limit theorem for the inferior arm at all.  The choice of a randomization schedule is, at the same time, a choice of how much can ever be learned about the inferior treatment.

Several directions remain open.  Natural next steps include finite-sample error bounds, delayed and batched responses, covariate-adjusted randomization, and Bayesian or regret-based analogues of the annealed schedules.  We leave these for future work.

\bibliographystyle{plainnat}
\bibliography{references}

\begin{thebibliography}{37}
\providecommand{\natexlab}[1]{#1}
\providecommand{\url}[1]{\texttt{#1}}
\expandafter\ifx\csname urlstyle\endcsname\relax
  \providecommand{\doi}[1]{doi: #1}\else
  \providecommand{\doi}{doi: \begingroup \urlstyle{rm}\Url}\fi

\bibitem[Anscombe(1952)]{Anscombe1952}
Francis~J. Anscombe.
\newblock Large-sample theory of sequential estimation.
\newblock \emph{Proceedings of the Cambridge Philosophical Society}, 48\penalty0 (4):\penalty0 600--607, 1952.
\newblock \url{https://doi.org/10.1017/S0305004100076386}.

\bibitem[Atkinson and Biswas(2014)]{AtkinsonBiswas2014}
Anthony~C. Atkinson and Atanu Biswas.
\newblock \emph{Randomised Response-Adaptive Designs in Clinical Trials}.
\newblock CRC Press, Boca Raton, FL, 2014.
\newblock \url{https://doi.org/10.1201/b16101}.

\bibitem[Auer et~al.(2002)Auer, Cesa-Bianchi, and Fischer]{AuerCesaBianchiFischer2002}
Peter Auer, Nicol\`o Cesa-Bianchi, and Paul Fischer.
\newblock Finite-time analysis of the multiarmed bandit problem.
\newblock \emph{Machine Learning}, 47:\penalty0 235--256, 2002.
\newblock \url{https://doi.org/10.1023/A:1013689704352}.

\bibitem[Baldi~Antognini and Giovagnoli(2015)]{BaldiAntogniniGiovagnoli2015}
Alessandro Baldi~Antognini and Alessandra Giovagnoli.
\newblock \emph{Adaptive Designs for Sequential Treatment Allocation}.
\newblock CRC Press, Boca Raton, FL, 2015.
\newblock \url{https://doi.org/10.1201/b18306}.

\bibitem[Bandyopadhyay and Biswas(2001)]{BandyopadhyayBiswas2001}
Uttam Bandyopadhyay and Atanu Biswas.
\newblock Adaptive designs for normal responses with prognostic factors.
\newblock \emph{Biometrika}, 88\penalty0 (2):\penalty0 409--419, 2001.
\newblock \url{https://doi.org/10.1093/biomet/88.2.409}.

\bibitem[Berger(2005)]{Berger2005}
Vance~W. Berger.
\newblock \emph{Selection Bias and Covariate Imbalances in Randomized Clinical Trials}.
\newblock Wiley, Chichester, 2005.
\newblock \url{https://doi.org/10.1002/0470863641}.

\bibitem[Berry and Fristedt(1985)]{BerryFristedt1985}
Donald~A. Berry and Bert Fristedt.
\newblock \emph{Bandit Problems: Sequential Allocation of Experiments}.
\newblock Chapman and Hall, London, 1985.

\bibitem[Bhandari et~al.(2007)Bhandari, Dutta, and Niyogi]{Bhandari2007}
Subir~Kumar Bhandari, R.~Dutta, and R.~G. Niyogi.
\newblock Study of optimal adaptive rule in testing problem.
\newblock In \emph{Advances in Multivariate Statistical Methods}. World Scientific, Singapore, 2007.

\bibitem[Bhandari et~al.(2009)Bhandari, De, Mandal, Pradhan, and Ghosh]{Bhandari2009}
Subir~Kumar Bhandari, S.~K. De, S.~Mandal, S.~Pradhan, and B.~Ghosh.
\newblock Study of optimal adaptive rule in testing composite hypothesis.
\newblock \emph{Sequential Analysis}, 28\penalty0 (3):\penalty0 394--405, 2009.
\newblock \url{https://doi.org/10.1080/07474940903041852}.

\bibitem[Biswas and Coad(2005)]{BiswasCoad2005}
Atanu Biswas and D.~Stephen Coad.
\newblock A general multi-treatment adaptive design for multivariate responses.
\newblock \emph{Sequential Analysis}, 24\penalty0 (2):\penalty0 139--158, 2005.
\newblock \url{https://doi.org/10.1081/SQA-200056191}.

\bibitem[Blackwell and Hodges(1957)]{BlackwellHodges1957}
David Blackwell and J.~L. Hodges.
\newblock Design for the control of selection bias.
\newblock \emph{The Annals of Mathematical Statistics}, 28\penalty0 (2):\penalty0 449--460, 1957.
\newblock \url{https://doi.org/10.1214/aoms/1177706973}.

\bibitem[Bubeck and Cesa-Bianchi(2012)]{BubeckCesaBianchi2012}
S\'ebastien Bubeck and Nicol\`o Cesa-Bianchi.
\newblock Regret analysis of stochastic and nonstochastic multi-armed bandit problems.
\newblock \emph{Foundations and Trends in Machine Learning}, 5\penalty0 (1):\penalty0 1--122, 2012.
\newblock \url{https://doi.org/10.1561/2200000024}.

\bibitem[Chow(1965)]{Chow1965}
Yuan~Shih Chow.
\newblock Local convergence of martingales and the law of large numbers.
\newblock \emph{The Annals of Mathematical Statistics}, 36\penalty0 (2):\penalty0 552--558, 1965.
\newblock \url{https://doi.org/10.1214/aoms/1177700166}.

\bibitem[Efron(1971)]{Efron1971}
Bradley Efron.
\newblock Forcing a sequential experiment to be balanced.
\newblock \emph{Biometrika}, 58\penalty0 (3):\penalty0 403--417, 1971.
\newblock \url{https://doi.org/10.1093/biomet/58.3.403}.

\bibitem[Eisele(1994)]{Eisele1994}
Jeffrey~R. Eisele.
\newblock The doubly adaptive biased coin design for sequential clinical trials.
\newblock \emph{Journal of Statistical Planning and Inference}, 38\penalty0 (2):\penalty0 249--261, 1994.
\newblock \url{https://doi.org/10.1016/0378-3758(94)90038-8}.

\bibitem[Freedman(1973)]{Freedman1973}
David~A. Freedman.
\newblock Another note on the borel-cantelli lemma and the strong law, with the poisson approximation as a byproduct.
\newblock \emph{The Annals of Probability}, 1\penalty0 (6):\penalty0 910--925, 1973.
\newblock \url{https://doi.org/10.1214/aop/1176996800}.

\bibitem[Friedman et~al.(2010)Friedman, Furberg, and DeMets]{FriedmanFurbergDeMets2010}
Lawrence~M. Friedman, Curt~D. Furberg, and David~L. DeMets.
\newblock \emph{Fundamentals of Clinical Trials}.
\newblock Springer, New York, 2010.
\newblock \url{https://doi.org/10.1007/978-1-4419-1586-3}.

\bibitem[Gut(2009)]{Gut2009}
Allan Gut.
\newblock \emph{Stopped Random Walks: Limit Theorems and Applications}.
\newblock Springer, New York, 2nd edition, 2009.
\newblock \url{https://doi.org/10.1007/978-0-387-87835-5}.

\bibitem[Hall and Heyde(1980)]{HallHeyde1980}
Peter Hall and Christopher~C. Heyde.
\newblock \emph{Martingale Limit Theory and Its Application}.
\newblock Academic Press, New York, 1980.

\bibitem[Hu and Rosenberger(2003)]{HuRosenberger2003}
Feifang Hu and William~F. Rosenberger.
\newblock Optimality, variability, power: Evaluating response-adaptive randomization procedures for treatment comparisons.
\newblock \emph{Journal of the American Statistical Association}, 98\penalty0 (463):\penalty0 671--678, 2003.
\newblock \url{https://doi.org/10.1198/016214503000000576}.

\bibitem[Hu and Rosenberger(2006)]{HuRosenberger2006}
Feifang Hu and William~F. Rosenberger.
\newblock \emph{The Theory of Response-Adaptive Randomization in Clinical Trials}.
\newblock Wiley, Hoboken, NJ, 2006.
\newblock \url{https://doi.org/10.1002/047005588X}.

\bibitem[Hu and Zhang(2004)]{HuZhang2004}
Feifang Hu and Li-Xin Zhang.
\newblock Asymptotic properties of doubly adaptive biased coin designs for multitreatment clinical trials.
\newblock \emph{The Annals of Statistics}, 32\penalty0 (1):\penalty0 268--301, 2004.
\newblock \url{https://doi.org/10.1214/aos/1079120137}.

\bibitem[Ivanova and Rosenberger(2000)]{IvanovaRosenberger2000}
Anastasia Ivanova and William~F. Rosenberger.
\newblock A comparison of urn designs for randomized clinical trials of {$k>2$} treatments.
\newblock \emph{Journal of Biopharmaceutical Statistics}, 10\penalty0 (1):\penalty0 93--107, 2000.
\newblock \url{https://doi.org/10.1081/BIP-100101016}.

\bibitem[Kundu et~al.(2025)Kundu, Jha, and Bhandari]{KunduJhaBhandari2025}
Sampurna Kundu, Jayant Jha, and Subir~Kumar Bhandari.
\newblock To study properties of a known procedure in adaptive sequential sampling design.
\newblock \emph{Sequential Analysis}, 44\penalty0 (4):\penalty0 451--463, 2025.
\newblock \url{https://doi.org/10.1080/07474946.2025.2533328}.

\bibitem[Lai and Robbins(1985)]{LaiRobbins1985}
Tze~Leung Lai and Herbert Robbins.
\newblock Asymptotically efficient adaptive allocation rules.
\newblock \emph{Advances in Applied Mathematics}, 6\penalty0 (1):\penalty0 4--22, 1985.
\newblock \url{https://doi.org/10.1016/0196-8858(85)90002-8}.

\bibitem[Lai and Wei(1982)]{LaiWei1982}
Tze~Leung Lai and Ching~Zong Wei.
\newblock Least squares estimates in stochastic regression models with applications to identification and control of dynamic systems.
\newblock \emph{The Annals of Statistics}, 10\penalty0 (1):\penalty0 154--166, 1982.
\newblock \url{https://doi.org/10.1214/aos/1176345697}.

\bibitem[Lattimore and Szepesv\'ari(2020)]{LattimoreSzepesvari2020}
Tor Lattimore and Csaba Szepesv\'ari.
\newblock \emph{Bandit Algorithms}.
\newblock Cambridge University Press, Cambridge, 2020.
\newblock \url{https://doi.org/10.1017/9781108571401}.

\bibitem[L\'evy(1937)]{Levy1937}
Paul L\'evy.
\newblock \emph{Th\'eorie de l'addition des variables al\'eatoires}.
\newblock Gauthier-Villars, Paris, 1937.
\newblock \url{https://doi.org/10.2307/3608623}.

\bibitem[Melfi and Page(2000)]{MelfiPage2000}
Vincent~F. Melfi and Connie Page.
\newblock Estimation after adaptive allocation.
\newblock \emph{Journal of Statistical Planning and Inference}, 87\penalty0 (2):\penalty0 353--363, 2000.
\newblock \url{https://doi.org/10.1016/S0378-3758(99)00198-6}.

\bibitem[Robbins(1952)]{Robbins1952}
Herbert Robbins.
\newblock Some aspects of the sequential design of experiments.
\newblock \emph{Bulletin of the American Mathematical Society}, 58\penalty0 (5):\penalty0 527--535, 1952.
\newblock \url{https://doi.org/10.1090/S0002-9904-1952-09620-8}.

\bibitem[Rosenberger et~al.(2001)Rosenberger, Stallard, Ivanova, Harper, and Ricks]{RosenbergerEtAl2001}
William~F. Rosenberger, Nigel Stallard, Anastasia Ivanova, Caroline~N. Harper, and Michelle~L. Ricks.
\newblock Optimal adaptive designs for binary response trials.
\newblock \emph{Biometrics}, 57\penalty0 (3):\penalty0 909--913, 2001.
\newblock \url{https://doi.org/10.1111/j.0006-341X.2001.00909.x}.

\bibitem[Sabatowski et~al.(2004)Sabatowski, Gálvez, Cherry, Jacquot, Vincent, Maisonobe, Versavel, and {The 1008-045 Study Group}]{SABATOWSKI2004}
Rainer Sabatowski, Rafael Gálvez, David~A Cherry, Florence Jacquot, Emmanuelle Vincent, Pascal Maisonobe, Mark Versavel, and {The 1008-045 Study Group}.
\newblock Pregabalin reduces pain and improves sleep and mood disturbances in patients with post-herpetic neuralgia: results of a randomised, placebo-controlled clinical trial.
\newblock \emph{Pain}, 109\penalty0 (1):\penalty0 26--35, 2004.
\newblock \url{https://doi.org/10.1016/j.pain.2004.01.001}.

\bibitem[Thompson(1933)]{Thompson1933}
William~R. Thompson.
\newblock On the likelihood that one unknown probability exceeds another in view of the evidence of two samples.
\newblock \emph{Biometrika}, 25\penalty0 (3--4):\penalty0 285--294, 1933.
\newblock \url{https://doi.org/10.1093/biomet/25.3-4.285}.

\bibitem[van~der Vaart(1998)]{vanderVaart1998}
Aad~W. van~der Vaart.
\newblock \emph{Asymptotic Statistics}.
\newblock Cambridge University Press, Cambridge, 1998.
\newblock \url{https://doi.org/10.1017/CBO9780511802256}.

\bibitem[Wei and Durham(1978)]{WeiDurham1978}
L.~J. Wei and Sue Durham.
\newblock The randomized play-the-winner rule in medical trials.
\newblock \emph{Journal of the American Statistical Association}, 73\penalty0 (364):\penalty0 840--843, 1978.
\newblock \url{https://doi.org/10.1080/01621459.1978.10480109}.

\bibitem[Williams(1991)]{Williams1991}
David Williams.
\newblock \emph{Probability with Martingales}.
\newblock Cambridge University Press, Cambridge, 1991.
\newblock \url{https://doi.org/10.1017/CBO9780511813658}.

\bibitem[Zelen(1969)]{Zelen1969}
Marvin Zelen.
\newblock Play the winner rule and the controlled clinical trial.
\newblock \emph{Journal of the American Statistical Association}, 64\penalty0 (325):\penalty0 131--146, 1969.
\newblock \url{https://doi.org/10.1080/01621459.1969.10500959}.

\end{thebibliography}

\end{document}